\documentclass[3p,11pt,authoryear]{elsarticle}
\pdfoutput=1
\usepackage[pagewise]{lineno}
\usepackage[hidelinks]{hyperref}

\usepackage{mathtools}

\DeclarePairedDelimiter\floor{\lfloor}{\rfloor}

\usepackage{longtable}
\usepackage{appendix}
\usepackage{bbm}
\usepackage{xcolor}
\usepackage{enumitem}
\usepackage[ruled]{algorithm2e}
\usepackage{adjustbox}
\usepackage{graphicx}
\usepackage{amsfonts}
\usepackage{amssymb}
\usepackage{amsmath}
\usepackage{multirow}
\usepackage{makecell}
\usepackage{rotating}
\usepackage{pdflscape}
\usepackage{arydshln}
\usepackage{etoolbox}
\usepackage{booktabs}
\usepackage{mathtools}
\usepackage{mathrsfs}
\usepackage{subcaption}
\usepackage{comment}
\usepackage{bm}
\usepackage{chngpage}
\usepackage{setspace}
\usepackage{url}
\usepackage{ltablex}
\usepackage{tabularx}

\usepackage{caption}
\usepackage{lscape}

\makeatletter 

\def\ps@pprintTitle{%
	\let\@oddhead\@empty
	\let\@evenhead\@empty
	\def\@oddfoot{\reset@font\hfil\thepage\hfil}
	\let\@evenfoot\@oddfoot
}
\makeatother

\makeatletter
\def\@author#1{\g@addto@macro\elsauthors{\normalsize%
		\def\baselinestretch{1}%
		\upshape\authorsep#1\unskip\textsuperscript{%
			\ifx\@fnmark\@empty\else\unskip\sep\@fnmark\let\sep=,\fi
			\ifx\@corref\@empty\else\unskip\sep\@corref\let\sep=,\fi
		}%
		\def\authorsep{\unskip,\space}%
		\global\let\@fnmark\@empty
		\global\let\@corref\@empty 
		\global\let\sep\@empty}%
	\@eadauthor={#1}
}
\makeatother

\begin{document}
\begin{frontmatter}

\title{Planning Charging Stations and Service Operations of Dockless Electric Micromobility Systems}




\author{Yining Liu}
\author{Yanfeng Ouyang\corref{correspondingauthor}}
\address{Department of Civil and Environmental Engineering, University of Illinois at Urbana-Champaign, Urbana, IL 61801, USA}
\cortext[correspondingauthor]{Corresponding author. Tel: +1 217 333 9858}
\ead{yfouyang@illinois.edu}

\begin{abstract}
Dockless electric micro-mobility services (e.g., shared e-scooters and e-bikes) have been increasingly popular in the recent decade, and a variety of charging technologies have emerged for these services. The use of charging stations, to/from which service vehicles are transported by the riders for charging, poses as a promising approach because it reduces the need for dedicated staff or contractors. However, unique challenges also arise, such as how to incentivize riders to drop off vehicles at stations and how to efficiently utilize the vehicles being charged at the stations. This paper focuses on dockless e-scooters as an example and develops a new spatial queuing network model to capture the steady-state scooter service cycles, battery consumption and charging processes, and the associated pricing and management mechanisms. Building upon this model, a system of closed-form equations is formulated and incorporated into a constrained nonlinear program to optimize the deployment of the service fleet, the design of charging stations (i.e., number, location, and capacity), user-based charging price promotions and priorities, and repositioning truck operations (i.e., headway and truck load). The proposed queuing network model is found to match very well with agent-based simulations. It is applied to a series of numerical experiments to draw insights into the optimal designs and the system performance. The numerical results reveal strong advantages of using charging stations for shared dockless electric micro-mobility services as compared to state-of-the-art alternatives. The proposed model can also be used to analyze other micromobility services and other charging approaches. 
\end{abstract}


\end{frontmatter}
\thispagestyle{empty}
\nolinenumbers

\section{Introduction}
The past decade has witnessed the emergence and
growth of shared electric micromobility services, in which e-bikes and e-scooters are used as a fun and convenient travel option for short-distance trips. These new services are expected to complement public transit systems, substitute automobile trips, alleviate congestion, and reduce emissions \citep{NABSA2023}. They have become very popular and their market sizes have been growing rapidly in the U.S. In 2018, shared e-scooters and shared e-bikes served 38.5 million trips and 6.5 million trips in the U.S., respectively, and these numbers quickly proliferated to 58.5 million and 20 million in 2022, respectively, despite the strong negative impacts of the COVID-19 pandemic \citep{NACTO2019,NACTO2022}. Shared electric micromobility services can be both docked (i.e., parking at docking stations only) and dockless (i.e., parking at any locations). The latter has been prevailing for shared e-scooter services in the U.S., while the former for the majority of shared e-bike services; however, dockless systems are also becoming popular for shared e-bikes -- the trip number increased by 73\% from 2.5 million in 2021 to 4.5 million in 2022. Overall, dockless systems serve 54\% of all U.S. micromobility trips in 2022 \citep{NACTO2022}.

The shared electric micromobility services, like any other shared mobility services, face many challenges. 
For example, for docked systems, the location (and capacity) of docking stations are critical planning decisions that must be optimized for maximal coverage 
\citep{frade2015bike}, minimal imbalance between pickups or drop-offs at each station \citep{liu2015station}, and/or minimal overall system-wide cost \citep{lin2011strategic}. In addition, 
a large number of shared electric bike/scooter fleet usually should be repositioned periodically 
to serve spatiotemporally heterogeneous demand. Studies have analyzed the spatiotemporal distribution patterns of demand for shared micromobility services \citep{zhu2020understanding,he2020dynamic}, determined optimal inventory level of shared fleet at each station of docked system \citep{raviv2013optimal,lu2016robust}, and proposed various strategies to reposition shared e-scooters or bikes, such as pricing incentives for riders to drop off vehicles at certain locations\citep{yun2022price,jin2023vehicle}, or repositioning bikes/scooters with dispatched trucks. On the latter, the assignment and routing of rebalancing trucks can be optimized with discrete modeling approaches \citep{chemla2013bike,o2015data,liu2016rebalancing,ghosh2017dynamic,ho2017hybrid,dell2018bike,warrington2019two}, with a discrete-continuous hybrid approach \citep{lei2018continuous,osorio2021optimal}, or with reinforcement learning \citep{losapio2021smart}. 
Existing literature also investigated other challenges of shared micromobility services, e.g., the competition between different service providers \citep{jiang2020optimal,jiang2022pricing}, the clutter of dockless shared fleet \citep{carrese2021beautiful}, and the integration with public transit \citep{wu2020optimal}.

The fact that shared electric micromobility vehicles must be recharged periodically introduces additional complexities. The battery level of a shared e-scooter or e-bike, also known as State-of-Charge (SoC), drops during its use. 
The charging process usually takes a non-trivial duration of time, which interrupts the vehicle service operations. It is hence crucial to determine the best approach to charge the shared electric fleet and optimize the charging operations. From a modeling point of view, vehicles at different SoCs are of different types (because they have different capabilities of providing future service) that should be treated differently in operation. 

There has been a large school of literature on the design of charging stations for privately owned electric vehicles, e.g., optimizing locations of charging stations along a corridor \citep{ghamami2016general} or on a road network \citep{kavianipour2021electric, he2018optimal}, considering the interactions with power systems \citep{lei2021system}, deploying multiple types of charging stations \citep{liu2017locating}, adopting dynamic charging infrastructure (e.g., charging lanes that charge vehicles while in motion) \citep{chen2017deployment}, and determining the density of charging stations in a uniform region via steady-state analysis \citep{varma2023electric}. 

More recently, researchers have also investigated the design of charging stations for docked micromobility services, e.g., shared e-bikes and shared electric vehicles. \citet{chen2020optimal} optimize the location and charging capacities of the dock stations of shared e-bikes as part of a discrete network in a bi-level model, where the lower level determines the riders' route choice based on the generalized cost (including the travel time, charging-related waiting time, and service fees). \citet{guo2022rebalancing} propose a Basket-Chandy-Muntz-Palacios (BCMP) queuing network to model the operation of station-based shared electric vehicles with charging delays, and optimize a variety of rebalancing strategies, including active rebalancing (by incentivizing riders to choose different destination stations) and passive rebalancing (by repositioning idle vehicles without serving a trip). \citet{yang2022user} explore another possibility in which dedicated trucks visit depleted shared e-bikes at dock stations and swap their batteries. To enhance the efficiency of battery swapping operation, incentive mechanisms are developed to encourage riders to drop off low-battery e-bikes at a subset of destination stations. 

Charging strategies have also been studied and implemented for dockless shared electric micromobility services. Companies, such as Bird and Lime, pay individual contractors (also known as ``chargers" for Bird and ``juicers" for Lime) a fee for rebalancing and charging an e-bike or an e-scooter. \citet{masoud2019optimal} optimize the assignment of e-scooters to individual contractors that minimizes their average travel distance via a mixed integer linear programming model. Another strategy is to hire technicians or ``swappers" to swap e-scooter/e-bike batteries in the field (e.g., adopted by Veoride). \citet{pender2020stochastic} build a queuing network model to analyze the operations of shared dockless e-scooters with swappable batteries under stochastic rider demand, trip durations, battery consumption processes, and swapper arrivals. A mean-field limit theorem and a functional central limit theorem are used to derive the percentage of e-scooters at each SoC and the number of swappers needed to guarantee any level of service. \citet{zhao2024planning} formulate a location-routing problem that optimizes the location of battery-swapping cabinets (that store depleted batteries and provide fully charged batteries) and the routing for technicians to visit depleted shared e-bikes. The concept of on-board charging emerged in recent years, which allows e-scooters/e-bikes to be charged on trucks during rebalancing operations. For example, \citet{osorio2021optimal} formulate a discrete vehicle routing problem to optimize the truck routes at a zone-level, and use a continuum approximation model to capture the routing cost for collecting and distributing e-scooters within each zone. SoC transitions of the e-scooters being charged on trucks are captured in both zone-level and local routes. 

More recently, the concept of planning charging stations for dockless shared e-scooters and e-bikes has been implemented in several major U.S. cities, e.g., New York and Atlanta \citep{stationAtlanta,stationNYC}. The charging stations can not only reduce the labor cost associated with independent contractors or dedicated crews, but also reshape the operations of dockless shared e-scooters and e-bikes by marrying the dockless system with the docked system --- now riders can pick up and drop off e-scooters/e-bikes at both random locations and charging stations, and as a result, are likely to have better access to service vehicles. These stations serve as both ``sink points" for the low-SoC e-scooter/e-bike, and supply points for the high-SoC ones, and hence alleviates the need to reposition these vehicles for charging. These benefits jointly lead to a higher system efficiency. However, on the other hand, such ``hybrid" systems also raise unique questions, such as (i) what are the probabilities for riders to select service vehicles of different SoCs among those at random locations or at the charging stations; (ii) how the riders would decide on the vehicle drop-off locations between their destinations or nearby stations; (iii) how the platform could incentivize riders to drop off vehicles (of different SoCs) at the stations for charging; (iv) how the platform could efficiently utilize e-scooters/e-bikes (of different SoCs) charged at stations; (v) how to plan the number, location, and capacity of these charging stations; and (vi) how can a platform conduct repositioning and recharging operations in case some depleted e-scooters/e-bikes are parked at the random locations. As such, planning charging stations for dockless e-scooters and e-bikes is not a trivial extension of the previous studies that solely focus on dockless systems and those on docked systems. The closest literature we can find, \citet{fathabad2022data}, adopts a discrete modeling approach to optimize the fleet size of e-scooters and the location and capacity of e-scooter charging stations, considering that e-scooters are rebalanced by either dedicated trucks or independent contractors to and from the charging stations. The dockless nature of shared e-scooter/e-bike services, however, is not explicitly captured. 


In light of the gaps in the literature, this paper proposes a new queuing network model for dockless shared electric micro-mobility services, and builds upon it to optimize the design of charging stations and service operations. To stay focused, we will take the shared e-scooter service as an example, but the models can be directly applied to e-bikes and other types of services. 

The proposed queuing network model is an extension of the modeling framework in \citet{Daganzo2019Paper}, which was originally used for analyzing shared mobility services (e.g., taxi, ride-pooling, dial-a-ride) in a homogeneous region, but later extended to services with maximum wait and detour guarantees \citep{Daganzo2020,ouyang2021performance}, integration with public transit \citep{liu2021mobility}, and services to spatially heterogeneous demand \citep{liu2023planning}. The model in this paper expands the queuing network state definition to capture the e-scooter SoC and incorporates transitions due to battery consumption and recharging. We allow the riders to pick up and drop off e-scooters of preferred SoC at both fixed locations (i.e., charging stations) and random locations, and allow the platform to use SoC-dependent promotions to incentivize riders to drop off depleted e-scooters at charging stations (i.e., user-based charging). 
The proposed system also uses a central charging depot, as a backup to the stations, to charge depleted e-scooters that are collected by trucks from random locations (i.e., truck-based charging). The system performance is evaluated by formulating a system of equations that includes a probabilistic description of scooter-rider matching at random locations and stations, expected e-scooter state durations, and state transition flow conservation. A constrained nonlinear program is built upon the system of equations to minimize the system-wide cost by optimizing a series of decisions, including the density and capacity of charging stations, the deployment of e-scooters at random locations and at stations, the SoC-based promotions for user-based charging, the truck repositioning headway and load, and the priority rule for e-scooters of different SoCs at stations to be available to riders.

The proposed queuing network model is verified with agent-based simulations under a range of parameter settings. Then, the performance of the proposed model is compared with two benchmark cases, a walk-only system (no e-scooters) and a depot-only system (in which all charging activities occur at the central depot), under a variety of hypothetical scenarios. The results demonstrate significant benefits of the proposed station-based system (especially when the platform puts high priorities for high-SoC e-scooters at stations to be used by riders) over the two benchmark cases for all investigated scenarios. 
As such, the proposed system, with properly planned service strategies and optimally deployed infrastructures, holds the promise to enhance operational efficiency of shared dockless electric micromobility services. More importantly, the proposed model can be easily modified to study other charging options (e.g., distributed charging hubs with no rider access, stations for swapping batteries, and hiring independent contractors), and further enable a comprehensive analysis of shared dockless electric micromobility services.

The remainder of the paper is organized as follows. Section \ref{scooter_sec:model} presents the proposed queuing network model for shared e-scooter services with charging stations, 
and the constrained nonlinear programming model. Section \ref{scooter_sec:numerical studies} first introduces the agent-based simulations, and then presents the numerical experiments as well as insights. Finally, Section \ref{scooter_sec:conclusion} concludes this study and discusses future research directions.

\section{Model}\label{scooter_sec:model}
In a square study region with side length $\Phi$ distance units (du), a transportation network company (TNC) deploys identical e-scooters 
to serve the travel demand. The city streets are assumed to form a very dense grid network, allowing e-scooter riders to travel along N-S and E-W directions and make a turn at any location; thus, the travel distance is measured by the rectilinear metric. The prorated cost of an e-scooter (e.g., for procurement and maintenance) per time unit (tu) is denoted by $\gamma$ [\$/tu]. Each e-scooter travels at a constant speed $v_s$ [du/tu], powered by identical batteries. The battery capacity of an e-scooter is discretized into a set of levels, also known as State-of-Charge (SoC), denoted by $\mathcal{B} = \{0, 1, 2, ..., B\}$. Without loss of generality, we assume each battery level allows an e-scooter to travel exactly one du, i.e., an e-scooter at SoC $b\in\mathcal{B}$ can travel up to $b$ du. For ease of explanation, we refer to an e-scooter at SoC $b$ as a type-$b$ e-scooter. 

The travel demand for e-scooters originates homogeneously in the study region at rate $\lambda$ [trips/tu-du$^2$], and all trips are local; i.e., no trips have a travel distance longer than a limit $L_{\max}$ [du], and $L_{\max}\ll\Phi$. The trip destination for each origin is uniformly distributed over the neighborhood within that range. 
As such, if we ignore the impacts of region boundary, the trip destinations are approximately uniformly distributed in the study region with rate $\lambda$. 
We further assume that $L_{\max}\leq B$ so that all trips can be served by fully charged e-scooters, and the riders do not change e-scooters in the middle of a trip. 

\begin{figure}
    \centering
    \includegraphics[width =0.8\textwidth]{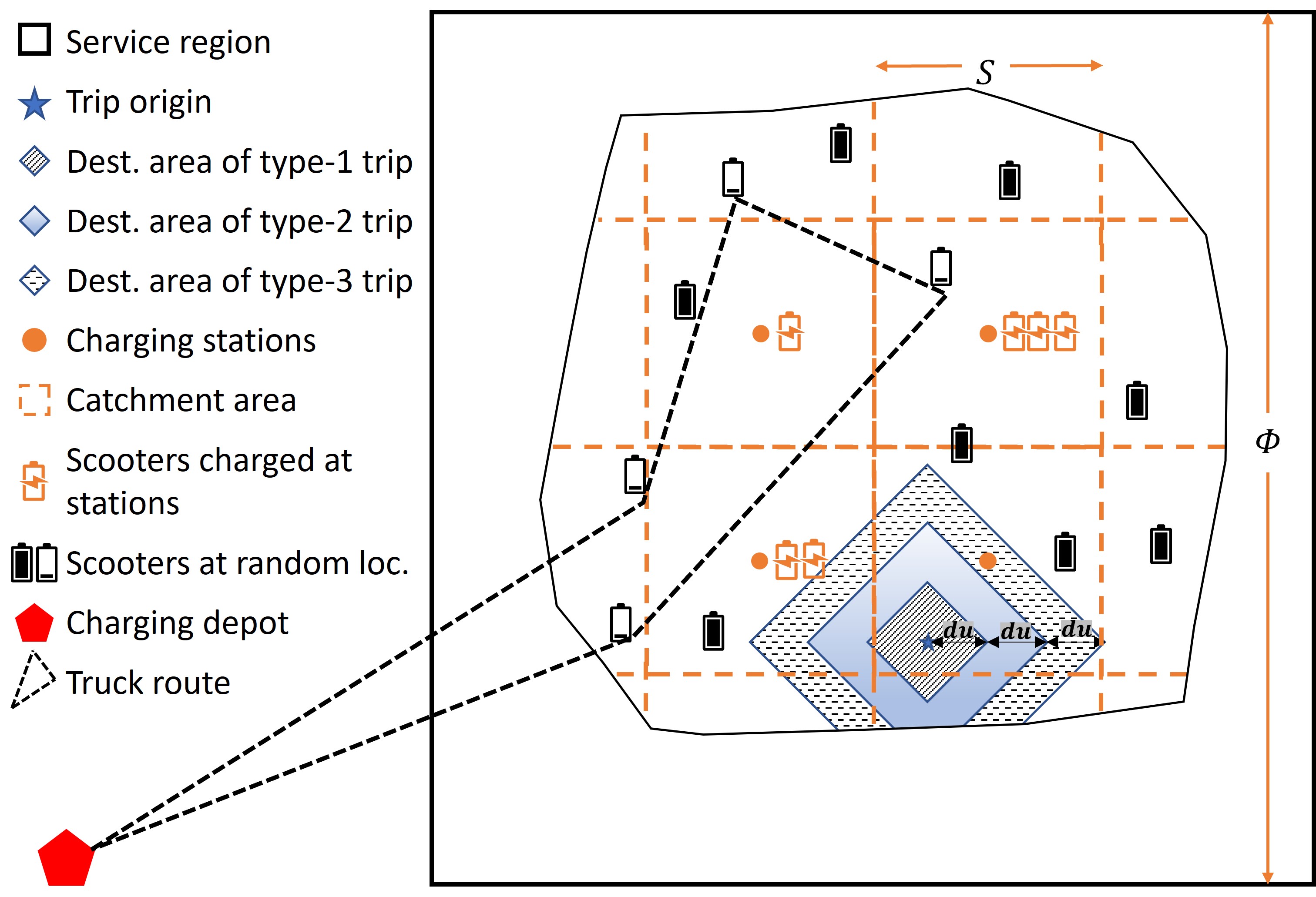}
    \caption{Illustration of shared e-scooter service with charging stations.}
    \label{scooter_fig:region illustration}
\end{figure}

We consider a trip with length $l$ du consumes $\hat b$ battery levels if $\hat{b}-1 < l \leq \hat{b}$, and this trip is called a type-$\hat b$ trip. It is easy to see that a trip consumes at least 1 battery level and at most $L_{\max}$ battery levels. Figure \ref{scooter_fig:region illustration} illustrates the destination areas of three types of trips from a single origin when $L_{\max} = 3$ du. For type-$\hat b$ riders, the trip rate $\lambda_{\hat{b}}$ [trips/tu-du$^2$] and the average trip length $L_{\hat{b}}$ [du] 
across all destinations can be calculated based on simple geometry, as follows,
\begin{equation}\label{scooter_eq: demand and trip length}
    \lambda_{\hat{b}} = \lambda \frac{(2\hat{b}-1)}{L^2_{\max}},\quad 
    L_{\hat{b}} = \frac{2}{3}\cdot\frac{3\hat{b}^2 - 3\hat{b} + 1}{2\hat{b} - 1},\quad \forall \hat{b}\in[1, 2, ..., L_{\max}].
\end{equation}
All riders have the same value-of-time $\beta$ [\$/tu] for both walking and riding, and they walk 
at a constant speed $v_w$ [du/tu]. The time or cost for a rider to unlock or return an e-scooter is assumed to be negligible. All key notation is summarized in \ref{scooter_sec:appendix notation}.

The TNC installs identical charging stations uniformly in the service region to form a grid with square catchment areas, and the distance between two adjacent stations is $S$ [du]; see Figure \ref{scooter_fig:region illustration} for an illustration. Their catchment areas equally partition the service area, such that $K\equiv\Phi/S$ is an integer, 
and there are $K^2$ 
charging stations in the study region. Every station is equipped with $Q$ chargers 
to provide simultaneous charging. The investments for each charging station and each charger, amortized to each operational time unit, are denoted by $\omega_1$ [\$/tu] and $\omega_2$ [\$/tu], respectively. 

The operating strategy is explained as follows. Some of the idle e-scooters are uniformly scattered in the service region, while others are randomly distributed at the charging stations. When a type-$\hat{b}$ rider searches for service, 
he/she always chooses a nearest suitable e-scooter (i.e., with sufficient battery) either at a random location or at the closest charging station.\footnote{We ignore the rare case in which the nearest suitable e-scooter is located at a farther charging station. This is reasonable when the number of idle e-scooters at random locations is relatively large.} The platform reserves this e-scooter for the rider upon request. 
In case multiple types of e-scooters are suitable and equally distant from this rider (e.g., all parked at a closest station), the platform imposes a priority rule for different types of e-scooters to be available to riders (e.g., by sorting the display sequence of, or even displaying certain suitable e-scooters as ``booked" in the app), to improve operation efficiency of charging stations. For example, the platform may prioritize e-scooters of higher SoC to be picked up by riders (so that they do not clog the charging stations) and prevent low-SoC ones from being picked up too frequently before adequate charging.\footnote{This priority rule can also be considered an input into the model to reflect the riders' preference for different types of e-scooters, which can be obtained from historical data.} 
Let $\theta_{\hat{b}, b}$ be a priority weight for a type-$b$ e-scooter at a station to be available to a type-$\hat{b}$ rider, $\forall \hat{b}\in[1, 2, ..., L_{\max}],\ b \in \mathcal{B}$, which satisfies: 
\begin{align}
    \theta_{\hat{b}, b} = 0,\ \forall b < \hat{b};\quad\quad  \theta_{\hat{b}, b} \geq 0,\ \forall b \geq \hat{b};\quad\quad \sum_{b \in \mathcal{B}} \theta_{\hat{b}, b} > 0,\ \forall \hat{b}.
\end{align}
The above formulas respectively define the priorities for unsuitable e-scooters and suitable e-scooters, and force the platform to provide at least one type of suitable e-scooters for any type of trips.
We assume the probability for a type-$\hat{b}$ rider to select a type-$b$ e-scooter out of equal-distance options is proportional to 
the number of type-$b$ e-scooters and the weight $\theta_{\hat{b}, b}$. 
For example, a rider is indifferent to any suitable e-scooters when $\theta_{\hat{b}, b} = 1,\ \forall b\geq \hat{b}$. 

After submitting the request, the type-$\hat{b}$ rider walks to the selected e-scooter of type $b$. The platform checks if the station closest to the rider's destination has unused or unbooked chargers. If so, the platform sends a promotion message (via the app) to encourage this rider to drop off the e-scooter at this charging station (i.e., user-based charging). The value of the promotion, denoted by $\pi_{b - \hat{b}}$ [\$/trip], depends on the estimated post-trip SoC of the e-scooter (i.e., $b - \hat{b}$). We assume that (i) the platform only incentivizes a rider to drop off an e-scooter at the station nearest to his/her destination; (ii) a rider would always accept the offer if the promotion value exceeds the disutility of walking from the charging station to his/her final destination; and (iii) if a rider accepts the offer, the platform will reserve a charging spot at the station. Otherwise, if there is no promotion or the rider rejects the offer, this rider will drop off the e-scooter at his/her destination. When a type-$b$ e-scooter is dropped off at a station, it is charged immediately, and it takes $\tau_{b}$ [tu] to reach SoC $b + 1$. If this e-scooter is booked by a new rider within this time, its charging process is considered interrupted, and its SoC will remain at $b$. 

Since the charging stations have limited capacities, the platform also uses a central charging depot, located $l_f$ [du] distance away from the service region center and with an infinite capacity, as the backup. Sufficiently large trucks are dispatched every $H$ [tu] headway, and each truck carries $R$ count of type-$B$ e-scooters (i.e., only fully-charged ones) upon leaving the depot, distributes them to random locations, and in the meantime collects randomly located type-$0$ e-scooters (i.e., only depleted ones) back to the depot for charging (i.e., truck-based charging). 
These trucks do not visit charging stations to pick up or drop off e-scooters, nor do they relocate e-scooters between two random locations.\footnote{We assume the charging stations to be only accessible by riders, to explore the advantage of the user-based charging option. This is not necessary --- the proposed model can be easily extended to allow trucks to relocate any idle e-scooters at any SoC level. In addition, we ignore the need of relocating e-scooters between random locations, because the travel demand is assumed to be spatially homogeneous. 
} 
Trucks travel at a constant speed $v_t$ [du/tu], and the operation cost of a truck per distance unit is denoted by $\kappa$ [\$/du]. 

As such, the platform integrates both user-based and truck-based charging operations, aiming to minimize the system-wide cost (i.e., agency cost and rider travel cost) by optimizing a set of planning and operational decisions, including the number of deployed e-scooters, the number and capacity of charging stations (i.e., $K^2$ and $Q$), the charging incentive $\{\pi_b\}$, the headway and initial load of trucks (i.e., $H$ and $R$), and the priority weights $\{\theta_{\hat{b}, b}\}$.

\subsection{Queuing Network Model}
\begin{figure}[t]
    \centering
    \includegraphics[width = \textwidth]{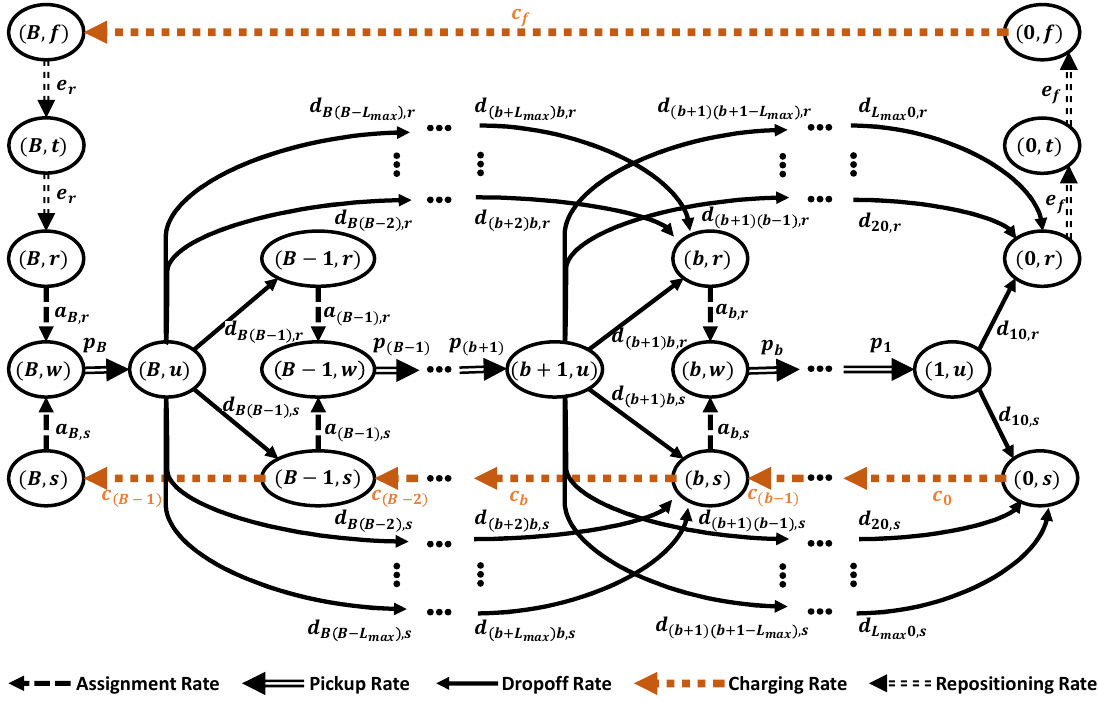}
    \caption{State Transition Network for Shared E-Scooter Service with Charging Stations}
    \label{scooter_fig:queuing network model}
\end{figure}

We develop a variant of the aspatial queuing network model proposed by \citet{Daganzo2019Paper} to capture the rider usage and the battery charging (at either the stations or the depot) of e-scooters in the steady state. Figure \ref{scooter_fig:queuing network model} gives a schematic illustration.

A possible state of the e-scooters is denoted as $(b, i)$, where $b\in\mathcal{B}$ represents the SoC and $i$ represents the service status. There are six possible service statuses: $i = w$ if the e-scooter is booked and waiting for pickup (i.e., the rider is walking toward the e-scooter), $i = u$ if the e-scooter is being used (i.e., it is carrying the rider toward the destination), or $i = r,\ s,\ t,\ f$ if the e-scooter is idle at a random location, being charged at a charging station, being transported by a truck, or being charged at the depot, respectively. All e-scooter states are represented by the circles in Figure \ref{scooter_fig:queuing network model}. The number of e-scooters in state $(b, i)$ is denoted by $n_{bi}$.

The state of an e-scooter changes when it is booked, picked up, dropped off, recharged, or repositioned. The state transitions are represented by the arrows in Figure \ref{scooter_fig:queuing network model}. We define $a_{b\hat{b},r}$ [trip/tu] and $a_{b\hat{b},s}$ [trip/tu] as the rates at which type-$b$ e-scooters at random locations and at stations are booked by type-$\hat{b}$ riders, respectively; the corresponding total booking rates of type-$b$ e-scooters (by all rider types) are $a_{b,r}=\sum_{\hat{b}\le b} a_{b\hat{b},r}$ and $a_{b,s}=\sum_{\hat{b}\le b} a_{b\hat{b},s}$, respectively, as illustrated by the dashed arrows; $p_b$ [trip/tu] the rate for type-$b$ e-scooters to be picked up, as illustrated by the double arrows; $d_{b(b-\hat{b}), s}$ [trip/tu] and $d_{b(b - \hat{b}), r}$ [trip/tu] the rates for type-$b$ e-scooters to be dropped off after serving type-$\hat{b}$ trips, at charging stations and at random locations, respectively, as illustrated by the solid arrows; $c_b$ [trip/tu] the rate for type-$b$ e-scooters to be charged into type-$(b+1)$ at charging stations, and $c_f$ [trip/tu] the rate for type-$0$ e-scooters to be charged into type-$B$ at the depot, as illustrated by the dotted arrow; and $e_f$ [trip/tu] and $e_r$ [trip/tu] the rates to reposition e-scooters from random locations to the depot, and from the depot to random locations, respectively, as illustrated by the double dashed arrows. 

Figure \ref{scooter_fig:queuing network model} illustrates all possible state transitions of shared e-scooters. For example, let's start with a type-$b$ e-scooter at a random location. If it is booked by a new type-$\hat{b}$ rider, it changes the service status from $i = r$ to booked $i = w$, i.e., its state $(b, r)$ is transitioned to state $(b, w)$ with rate $a_{b,r}$. Then, the rider walks to pick up the e-scooter. Upon pickup, the e-scooter updates the service status from booked $i=w$ to in-service $i = u$, i.e., it transitions to state $(b, u)$ with rate $p_b$. Then, the rider starts to ride the e-scooter. Based on the promotion $\pi_{\left(b-\hat{b}\right)}$, he/she may drop off the e-scooter either at its destination (a random location), which occurs at rate $d_{b(b-\hat{b}), r}$, or the station nearest to its destination, which occurs at rate $d_{b(b-\hat{b}), s}$. In the former case, upon drop-off, the e-scooter service status changes from in-service $i=u$ to idle at a random location $i=r$, and the e-scooter SoC drops based on the battery levels consumed by the trip $\hat{b}$, i.e., the e-scooter transitions to state $(b-\hat{b}, r)$. In the latter case, the e-scooter transitions to state $(b-\hat{b}, s)$, and starts to be charged toward SoC $b - \hat{b} + 1$. During the charging process, this e-scooter may be booked by a new rider again and transitions immediately to state $(b-\hat{b},w)$,\footnote{For simplicity, we ignore the possibility of continuing charging during the relatively short wait time between booking and pickup.} which occurs at rate $a_{(b-\hat{b}),s}$; otherwise, it transitions to state $(b-\hat{b}+1,s)$, which occurs at rate $c_{(b-\hat{b})}$. The type-$0$ e-scooters accumulated at random locations, in state $(0, r)$, are collected by trucks at flow rate $e_f$. Their state changes to $(0, t)$ upon pickup by the trucks, and $(0, f)$ upon arrival at the depot. They are charged at the depot until reaching full SoC $B$, transitioning to state $(B, f)$ at rate $c_f$, and waiting for the next truck dispatch. Then, these e-scooters are distributed by trucks toward random locations, in state $(B, t)$, and the transition rate is $e_r$.

This queuing network model has two sets of variables: flow rates $\{a_{b\hat{b},r}\}$, $\{a_{b\hat{b},s}\}$, $\{p_b\}$, $\{d_{b(b-\hat{b}), r}\}$, $\{d_{b(b-\hat{b}), s}\}$, $\{c_b\}$, $c_f$, $\{e_i\}$, and the number of e-scooters at each state $\{n_{bi}\}$. We next quantify their relationship by building a system of equations based on the probabilities of riders' e-scooter choice and drop-off location decisions, the expected durations of e-scooter waiting, service, and charging processes, and the flow conservation. 

\subsubsection{Booking, Pickup, Drop-off, and Charging Rates}\label{scooter_sec:booking probability}
A type-$\hat{b}$ rider's suitable idle e-scooters (whose SoC $b \geq \hat{b}$) can either be located at charging stations or at random locations. Their quantities are 
\begin{equation}
    N_{\hat{b}s} = \sum_{b : \theta_{\hat{b}, b} > 0} n_{bs}
    ,\quad N_{\hat{b}r} = \sum^{B}_{b = \hat{b}} n_{br}, \quad \forall \hat{b}\in[1, 2, ..., L_{\max}].
\end{equation}
The total number of e-scooters at all stations is $N_s = \sum_{b\in\mathcal{B}} n_{bs}$.

For a type-$\hat{b}$ rider to choose a suitable e-scooter at the nearest station, the following two independent conditions must hold at the same time: (i) some suitable e-scooters are available at this station, referred to as condition $C1$, which occurs with probability $P_{C1, \hat{b}}$; and (ii) no suitable e-scooters are available at random locations closer than this station, referred to as condition $C2$, which occurs with probability $P_{C2, \hat{b}}$. We can derive $P_{C1, \hat{b}}$ from its complementary probability (i.e., for all e-scooters at the nearest station to be unsuitable). Among all e-scooters that are at the $K^2$ stations, at most $Q$ e-scooters can be at this nearest station, as well as at each of the other $K^2-1$ stations; thus, the number of e-scooters at this station can neither exceed $Q$ nor be lower than $\max\{0,\ N_s-(K^2-1)Q\}$. 
When there are $q$ e-scooters at this station, there are no suitable e-scooters for this rider if and only if all these $q$ e-scooters are among the $N_s - N_{\hat{b}s}$ unsuitable e-scooters, and all remaining e-scooters are distributed at other stations; the probability for this to occur is approximately $\binom{N_s - N_{\hat{b}s}}{q}\left(\frac{1}{K^2}\right)^q \left(1 - \frac{1}{K^2}\right)^{N_s - q}$.\footnote{The calculation here approximates the capacity constraint of other $K^2-1$ stations based on their total capacities, while the capacity constraint of each station is not guaranteed. This approximation is reasonable especially when $N_s \ll K^2Q$, as verified by the agent-based simulations in Section \ref{scooter_sec:simulation}.} Therefore, 
\belowdisplayskip=20pt
\begin{equation}\label{scooter_eq:station booking condition 1}
\begin{aligned}
    P_{C1, \hat{b}} = 1 - \sum_{q = \max\{0,\ N_s-(K^2-1)Q\}}^{Q}\binom{N_s - N_{\hat{b}s}}{q}\left(\frac{1}{K^2}\right)^q \left(1 - \frac{1}{K^2}\right)^{N_s - q}.
\end{aligned}
\end{equation}

For condition $C2$, let random variable $l$ denote the distance from a type-$\hat{b}$ rider's origin to its nearest station. Conditional on $l$, a suitable e-scooter at a random location may be farther than this nearest station with a probability $1 - 2l^2/\Phi^2$. Therefore, 
none of the $N_{\hat{b}r}$ suitable e-scooters would be at a closer random location with probability $\left(1 - 2l^2/\Phi^2 \right)^{N_{\hat{b}r}}$.
Variable $l$ itself follows a triangular distribution, with limits $0$ and $S$, and mode $S/2$. Then, 
\begin{subequations}\label{scooter_eq:prob station closer by steps}
\renewcommand{\theequation}{\theparentequation.\arabic{equation}}
\begin{align}
    P_{C2, \hat{b}} &= 
    \int_{0}^{S/2}\left(1 - \frac{2l^2}{\Phi^2} \right)^{N_{\hat{b}r}}\frac{4l}{S^2}dl + \int_{S/2}^{S}\left(1 - \frac{2l^2}{\Phi^2} \right)^{N_{\hat{b}r}}\frac{4(S-l)}{S^2}dl \label{scooter_eq: prob station closer uncondition numerical}\\
    & = \frac{2K^2}{N_{\hat{b}r}+1}\left[1 - 2\left(1 - \frac{1}{2K^2}\right)^{N_{\hat{b}r}+1} + \left(1 - \frac{2}{K^2}\right)^{N_{\hat{b}r}+1}\right] + \sum^{N_{\hat{b}r}}_{i=0} \binom{N_{\hat{b}r}}{i}\frac{2 - 0.25^i}{i+0.5}\left(\frac{-2}{K^2}\right)^i.\label{scooter_eq: prob station closer uncondition}
\end{align}
\end{subequations}
At this station, the rider may book a specific suitable type-$b$ e-scooter based on the priority weights $\{\theta_{\hat{b}, b}\}$ and the numbers of e-scooters at charging stations $\{n_{bs}\}$, i.e., $\frac{\theta_{\hat{b}, b}\cdot n_{b s}}{\sum_{b'} \theta_{\hat{b}, b'}\cdot n_{b's}}$. Hence, the probability for a type-$\hat{b}$ rider to choose a suitable type-$b$ e-scooter at the nearest station is $P_{C1, \hat{b}}\cdot P_{C2, \hat{b}} \cdot \frac{\theta_{\hat{b}, b}\cdot n_{b s}}{\sum_{b'} \theta_{\hat{b}, b'}\cdot n_{b's}}$. 

On the other hand, this rider's nearest suitable e-scooter may be at a random location, and the probability of selecting a type-$b$ e-scooter is proportional to the relative density of type-$b$ e-scooters among all those suitable ones at random locations, i.e., $n_{br}/N_{\hat{b}r}$. Then, the rates for a type-$\hat{b}$ rider to book a type-$b$ e-scooter 
are
\belowdisplayskip=10pt
\begin{subequations}
\renewcommand{\theequation}{\theparentequation.\arabic{equation}}
\begin{align}
    &a_{b\hat{b}, s} = P_{C1, \hat{b}}P_{C2, \hat{b}}\cdot \frac{\theta_{\hat{b}, b}\cdot n_{bs}}{\sum_{b'} \theta_{\hat{b}, b'}\cdot n_{b's}}\cdot \lambda_{\hat{b}}\Phi^2,\\
    &\ a_{b\hat{b}, r} = \begin{cases}
        \left(1 - P_{C1, \hat{b}} P_{C2, \hat{b}}\right) \frac{n_{br}}{N_{\hat{b}r}}\cdot \lambda_{\hat{b}}\Phi^2, &\forall b \geq \hat{b}\\
        0, & \forall b < \hat{b}
    \end{cases},
\end{align}
\end{subequations}
and 
\begin{equation}
    a_{b, s} = \sum^{L_{\max}}_{\hat{b} = 1} a_{b\hat{b}, s},\quad a_{b, r} = \sum^{L_{\max}}_{\hat{b} = 1} a_{b\hat{b}, r},\quad\quad \forall b\in\mathcal{B}.
\end{equation}

\begin{figure}
     \centering
     \begin{subfigure}[b]{0.7\textwidth}
         \centering
         \includegraphics[width=\textwidth]{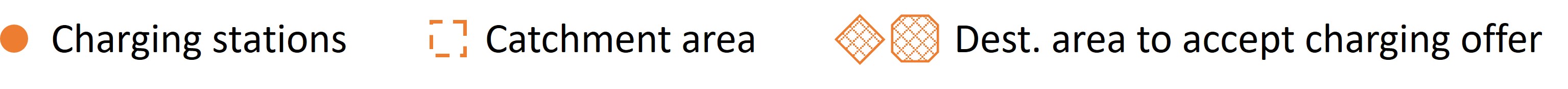}
     \end{subfigure}
     \hfill
     \begin{subfigure}[b]{0.4\textwidth}
         \centering
         \includegraphics[width=\textwidth]{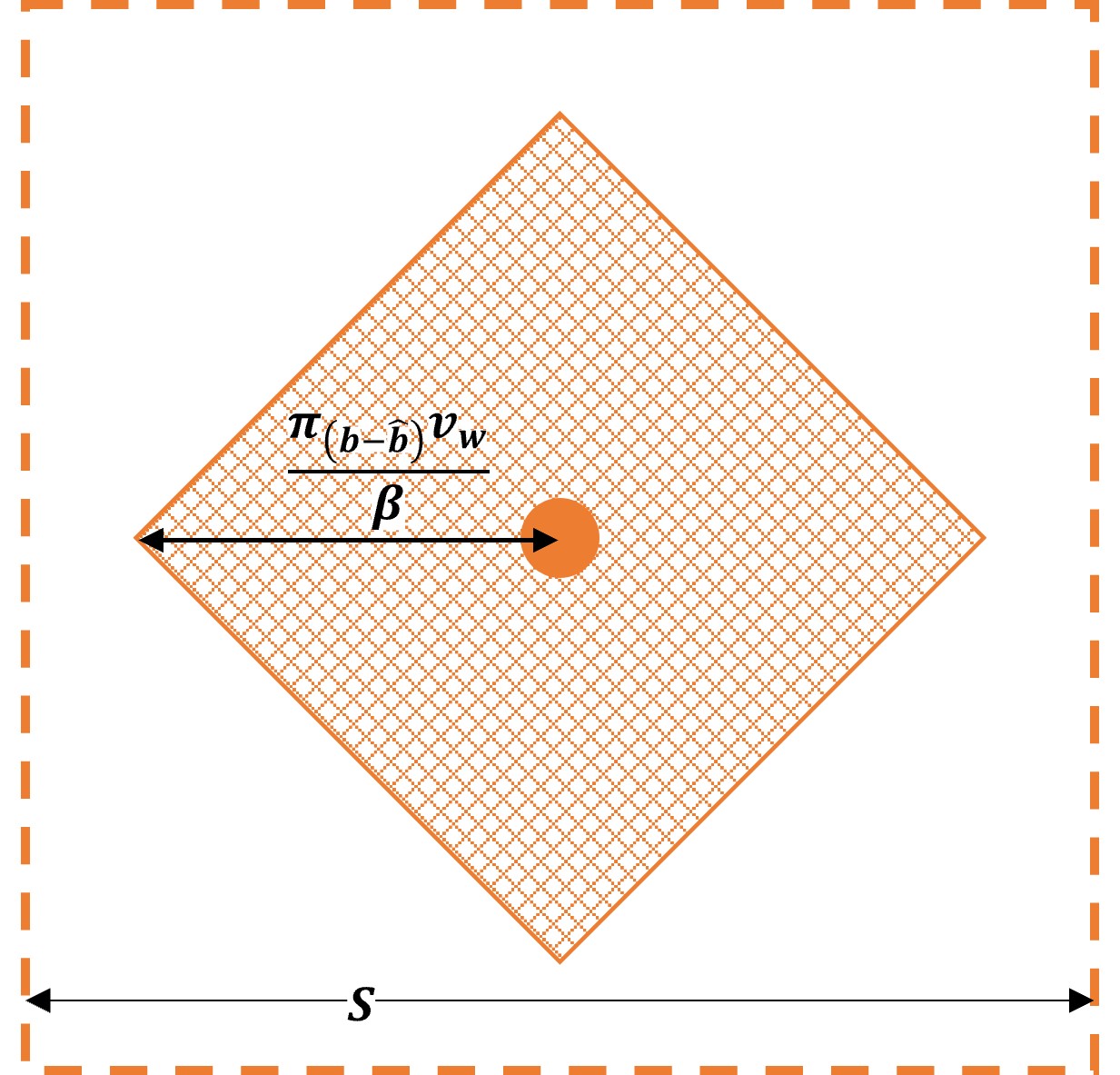}
         \caption{When $\pi_{(b-\hat{b})} \leq \frac{\beta S}{2v_w}$.}
         \label{scooter_fig:dest area pi 1}
     \end{subfigure}
     \begin{subfigure}[b]{0.4\textwidth}
         \centering
         \includegraphics[width=\textwidth]{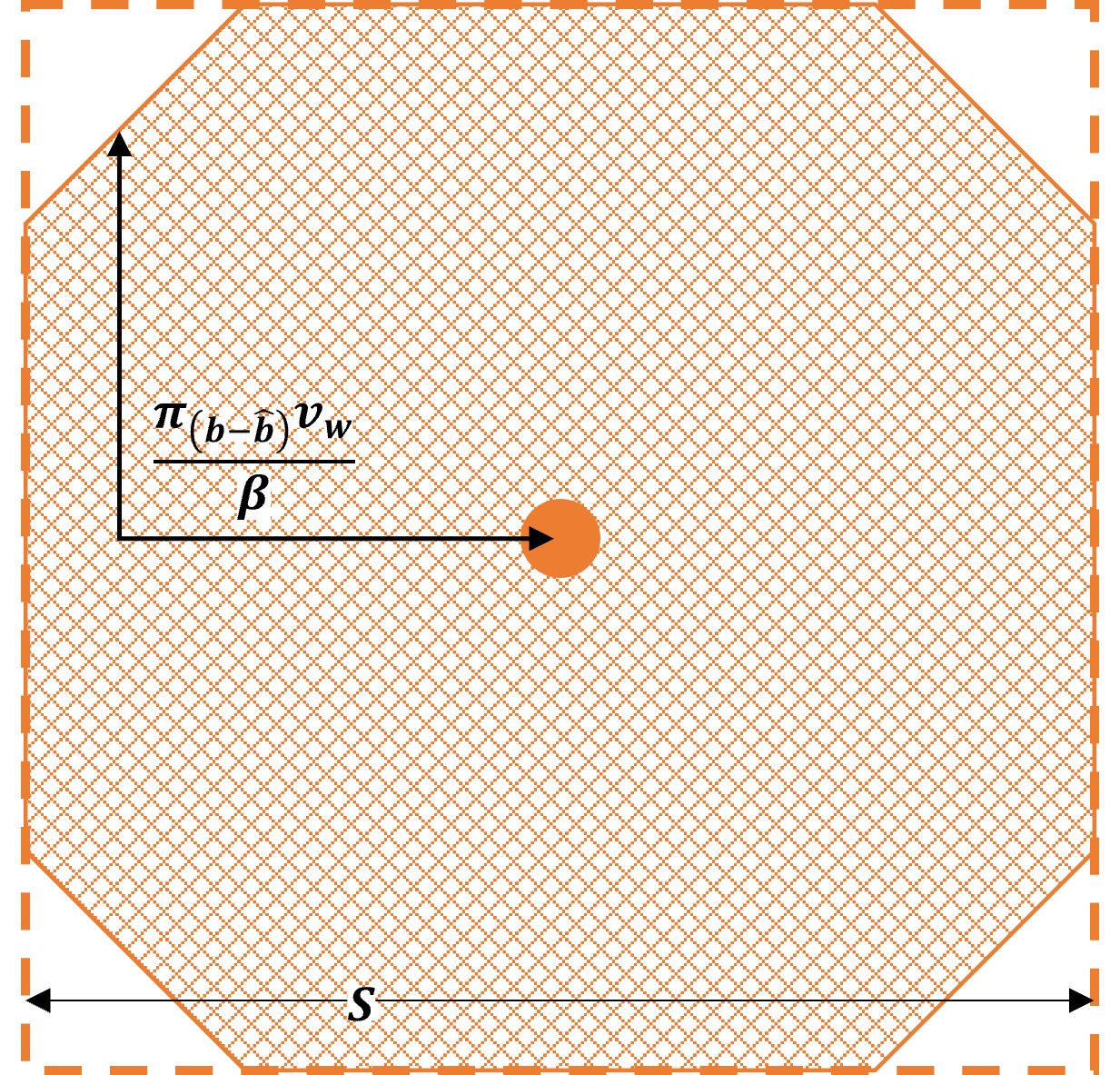}
         \caption{When $\frac{\beta S}{2v_w} < \pi_{\left(b-\hat{b}\right)} \leq \frac{\beta S}{v_w}$.}
         \label{scooter_fig:dest area pi 2}
     \end{subfigure}
    \caption{Illustration of the destination area for a type-$\hat{b}$ rider on a type-$b$ e-scooter to accept the charging incentive.}
    \label{scooter_fig:dest area pi}
\end{figure}

A rider drops off an e-scooter at the station nearest to its destination if and only if both of the following two independent conditions are satisfied: (i) this station has vacant chargers; and (ii) the promotion value exceeds the extra walking disutility. The probability of the first condition, 
denoted by $P_Q$, 
can be approximated with a binomial distribution, following the same idea of Eqn. \eqref{scooter_eq:station booking condition 1}, as follows:
\begin{equation}
    P_Q = \sum_{q = \max\{0,\ N_s-(K^2-1)Q\}}^{Q-1}\binom{N_s}{q}\left(\frac{1}{K^2}\right)^q \left(1 - \frac{1}{K^2}\right)^{N_s - q}.
\end{equation}
The second condition holds, 
for an e-scooter with post-trip SoC $b-\hat{b}$, 
as long as the distance between this rider's destination and the nearest station is less than $\pi_{\left(b-\hat{b}\right)} v_w/\beta$; see Figure \ref{scooter_fig:dest area pi} for an illustration.
\footnote{It is noted that under this strategy, no type-$(b-\hat{b})$ e-scooters would exist within the shaded areas, which violates the uniformity assumption for every type of idle e-scooters. 
This approximation is acceptable especially when $\pi_{b-\hat{b}}$ is either small (i.e., the shaded area in Figure \ref{scooter_fig:dest area pi} is small) or large such that $\pi_{\left(b-\hat{b}\right)} v_w/\beta \approx S$ (i.e., most of e-scooters are dropped off at stations and 
thus the distribution of state-$(b-\hat{b}, r)$ e-scooters does not significantly affect the system performance). 
This is verified by our numerical simulations in Section \ref{scooter_sec:simulation}.} Then, the probability for the charging promotion of a type-$b-\hat{b}$ e-scooter to exceed the rider's walking disutility, denoted by $P_{\pi,(b-\hat{b})}$, is, 
\begin{equation}
    P_{\pi,(b-\hat{b})} = \begin{cases}
        2\left(\frac{\pi_{\left(b-\hat{b}\right)} v_w}{\beta S}\right)^2,\quad\quad & \pi_{\left(b-\hat{b}\right)} \leq \frac{S\beta}{2v_w}\\
        1 - 2\left(1 - \frac{\pi_{\left(b-\hat{b}\right)} v_w}{\beta S}\right)^2,\quad\quad & \frac{\beta S}{2v_w} < \pi_{\left(b-\hat{b}\right)} \leq \frac{\beta S}{v_w}
    \end{cases}\,.
\end{equation}
Therefore, the drop-off rates for type-$(b-\hat{b})$ e-scooters at the stations and random locations are, respectively, the following:
\begin{equation}\label{scooter_eq: delivery rates}
    d_{b(b-\hat{b}),s} = P_Q P_{\pi,(b-\hat{b})} \cdot\left(a_{b\hat{b}, s} + a_{b\hat{b}, r}\right),\quad d_{b(b-\hat{b}),r} = a_{b\hat{b}, s} + a_{b\hat{b}, r} - d_{b(b-\hat{b}),s},\quad\quad \forall b\geq\hat{b}.
\end{equation}

When a type-$b$ e-scooter is being charged at a station, the probability of receiving no new service bookings during its charging time $\tau_b$ is $e^{-\frac{a_{b,s}}{n_{bs}}\cdot\tau_b}$, given Poisson arrivals of bookings at rate $a_{b,s}$ for all e-scooters in state $(b, s)$. 
In addition, type-$b$ e-scooters being charged at the stations come from (i) type-$b$ e-scooters that are dropped off
, at rate $\sum_{\hat{b}} d_{(b+\hat{b})b, s}$; and (ii) type-$(b-1)$ e-scooters that are charged to Soc $b$, at rate $c_{b-1}$. Hence, the rate at which type-$b$ e-scooters are charged to SoC $(b+1)$ is,
\begin{equation}
    \begin{aligned}
        c_{b} = \left[\displaystyle \sum^{\min\{L_{\max}, B - b\}}_{\hat{b} = 1}d_{(b + \hat{b})b,s} + c_{(b-1)}\right]\cdot e^{-\frac{a_{b,s}}{n_{bs}}\cdot \tau_b}, \quad & \forall b\in\{1, 2, ..., B-1\}.
    \end{aligned}
\end{equation}

\subsubsection{Flow Conservation}
In Figure \ref{scooter_fig:queuing network model}, the transition flows in the queuing network should be conserved at each state node, as described by the following system of equations.
\begin{subequations}\label{scooter_eq:flow conservation}
\renewcommand{\theequation}{\theparentequation.\arabic{equation}}
\begin{align}
&0 = c_{b} + a_{b,s} - \sum_{\hat{b} = 1}^{\min\{L_{\max}, B - b\}}d_{(b+\hat{b})b,s} - c_{(b-1)},\quad &&\forall b\in\{1, 2, ..., B-1\} \label{scooter_eq: flow conservation 1}\\
&0 = c_{0} - \sum_{\hat{b} = 1}^{L_{\max}}d_{\hat{b} 0,s} \label{scooter_eq: flow conservation 2}\\
&0 = a_{B,s} - c_{(B-1)} \label{scooter_eq: flow conservation 3}\\
&0 = a_{b,r} - \sum_{\hat{b} = 1}^{\min\{L_{\max}, B - b\}}d_{(b+\hat{b})b,r},\quad && \forall b\in\{1, 2, ..., B-1\} \\
&0 = e_f - \sum_{\hat{b} = 1}^{L_{\max}}d_{\hat{b} 0,r}\\
&0 = a_{B,r} - e_r\\
&0 = p_{b} - a_{b,s} - a_{b,r},\quad && \forall b\in\{1, 2, ..., B\}\\
&0 = \sum_{\hat{b} = 1}^{\min\{L_{\max}, b\}}\left(d_{b(b-\hat{b}),s} + d_{b(b-\hat{b}),r}\right) - p_{b},\quad && \forall b\in\{1, 2, ..., B\}\\
&c_f = e_f = e_r \label{scooter_eq: repositioning equation}.
\end{align}
\end{subequations}

\subsubsection{State Duration and Fleet Size}
Section \ref{scooter_sec:booking probability} shows that the flow rates $\{a_{b\hat{b},r}\}$, $\{a_{b\hat{b},s}\}$, $\{p_b\}$, $\{d_{b(b-\hat{b}), r}\}$, $\{d_{b(b-\hat{b}), s}\}$, $\{c_b\}$, $c_f$, and $\{e_i\}$ are dependent on only the numbers of idle e-scooters at stations $\{n_{bs}\}$ and at random locations $\{n_{br}\}_{b > 0}$, station number $K^2$, station capacity $Q$, promotions $\{\pi_b\}$, and priority weights $\{\theta_{\hat{b}, b}\}$. 
Once we know these rates, the number of e-scooters in all other states, i.e., $\{(b, w)\}$, $\{(b, u)\}$, $(0, r)$, $(0, t)$, $(0, f)$, $(B, t)$, and $(B, f)$, can be obtained based on the expected durations of these states per Little's Law \citep{Little1961}. This is shown next.

The expected pickup distance for a type-$\hat{b}$ rider, denoted by $l_{\hat{b},p}$, 
equals $0.63\Phi/\sqrt{N_{\hat{b}r}}$ \citep{Daganzo2019Paper} when there are no suitable e-scooters at the nearest station (i.e., condition $C1$ does not hold). 
Otherwise, when there are some suitable e-scooters at the nearest station (i.e., condition $C1$ holds), 
the pickup distance $l_{\hat{b},p}$ should not exceed the distance to this nearest station, $l$. 
It is straightforward to derive the 
CDF of $l_{\hat{b},p}$ conditioning on $l$ and the occurrence of $C1$, as follows: 
\begin{equation}\label{scooter_eq:lp cdf}
    1-F_{l_{\hat{b},p}}\left(x|l, C1\right) = \Pr\left\{l_{\hat{b},p} > x|l, C1\right\} = \begin{cases}
    \left(1-\frac{2x^2}{\Phi^2}\right)^{N_{\hat{b}r}}, & x \leq l \\
    0, & x > l
  \end{cases},
\end{equation}
and hence
\begin{subequations}\label{scooter_eq:lp by steps}
\renewcommand{\theequation}{\theparentequation.\arabic{equation}}
\small
\begin{align}
    \mathbb{E}\left[l_{\hat{b},p}|l, C1\right] &= \int_{0}^{l}x dF_{l_{\hat{b},p}}\left(x|l, C1\right) + l \left(1-\frac{2l^2}{\Phi^2}\right)^{N_{\hat{b}r}} = \int_0^l \left(1 - \frac{2x^2}{\Phi^2}\right)^{N_{\hat{b}r}}dx,\\
    \mathbb{E}\left[l_{\hat{b},p}|C1\right] &=\int_{0}^{\frac{S}{2}}\int_0^l \left(1 - \frac{2x^2}{\Phi^2}\right)^{N_{\hat{b}r}}\cdot \frac{4l}{S^2}dxdl + \int_{\frac{S}{2}}^{S}\int_0^l \left(1 - \frac{2x^2}{\Phi^2}\right)^{N_{\hat{b}r}}\cdot\frac{4(S-l)}{S^2}dxdl,  \label{scooter_eq:pickup dist numerical}\\
    &= \sum_{i = 0}^{N_{\hat{b}r}}{N_{\hat{b}r} \choose i}\frac{(-1)^i(2^{i+2} - 2^{-i})\cdot S}{(2i+1)(2i+2)(2i+3)\cdot K^{2i}}.
    \label{scooter_eq:pickup dist}
\end{align}
\end{subequations}
Then, the overall expected pickup distance for a type-$\hat{b}$ rider is
\begin{equation}\label{scooter_eq:lp final}
\small
    \mathbb{E}\left[l_{\hat{b},p}\right] = P_{C1, \hat{b}}\sum_{i = 0}^{N_{\hat{b}r}}{N_{\hat{b}r} \choose i}\frac{(-1)^i(2^{i+2} - 2^{-i})\cdot S}{(2i+1)(2i+2)(2i+3)\cdot K^{2i}}+ (1 - P_{C1, \hat{b}})\frac{0.63\Phi}{\sqrt{N_{\hat{b}r}}}.
\end{equation}
For an e-scooter in service state $(b, u)$, 
the expected service distance is approximately the same as this rider's trip length; i.e., $L_{\hat{b}}$ as in Eqn. \eqref{scooter_eq: demand and trip length} for a type-$\hat{b}$ rider.

Therefore, the numbers of e-scooters in states $(b, w)$ and $(b, u)$, based on the Little's Law (\citeyear{Little1961}), are as follows,
\begin{equation}
    n_{bw} = \sum_{\hat{b} = 1}^{\min\{L_{\max}, b\}} (a_{b\hat{b}, s} + a_{b\hat{b} ,r})\cdot \frac{E\left[l_{\hat{b},p}\right]}{v_w},\quad n_{bu} = \sum_{\hat{b} = 1}^{\min\{L_{\max}, b\}} \left[d_{b(b-\hat{b}), s} + d_{b(b-\hat{b}), r}\right]\cdot \frac{L_{\hat{b}}}{v_s}.
\end{equation}

The expected amounts of time e-scooters spend in states $(0, r)$, $(0, t)$, $(0, f)$, $(B, f)$ and $(B, t)$ are 
closely related to the repositioning headway $H$ and the length of truck routes. For every repositioning headway $H$, an average of $(H\cdot e_f)$ type-$0$ e-scooters and $(H\cdot e_r)$ type-$B$ e-scooters are transported into and out of the depot, respectively. Since each truck carries $R$ number of type-$B$ e-scooters out of the depot, the number of trucks needed 
is $m = H\cdot e_r/R= H\cdot e_f/R$. The optimal truck route length can be given by the asymptotic vehicle routing problem (VRP) 
formula \citep{beardwood_halton_hammersley_1959,ROBUST1990263} based on the depot distance $l_f$ and e-scooter densities: $l_e = 2ml_f + 0.95\sqrt{\Phi^2 (He_r + He_f)}$.
\footnote{For simplicity, we assume that $l_f \gg \Phi$, such that the optimal repositioning headway $H$ and the initial truck load $R$ are approximately the same for the entire service region.} The average route distance per truck per headway is $l_e/m$.

The expected waiting time for a randomly located type-$0$ e-scooter in state $(0, r)$ to be picked up by a truck 
includes the expected time to wait for the truck dispatch $H/2$ and for the truck to pick up and drop off other e-scooters before picking up this e-scooter $l_e/(2mv_t)$. 
The expected in-vehicle travel time is 
$l_e/(2mv_t)$. The charge duration from state-$(0, f)$ to state-$(B, f)$ is $\sum_b \tau_b$. 
Type-$B$ e-scooters at the charging facility 
first wait $H/2$ time on average for pickup by trucks, and then they 
spend $l_e/(2mv_t)$ time on the trucks before drop-off. 
Per Little's Law \citeyearpar{Little1961}, the following holds: 
\begin{align}\label{scooter_eq:fleet size}
n_{0r} = \frac{e_fH}{2} + \frac{e_fl_e}{2mv_t},\quad n_{0t} = \frac{e_fl_e}{2mv_t},\quad n_{0f} = \left(\sum^{B - 1}_{b=0} \tau_b\right)\cdot c_f,\quad
n_{Bf} = \frac{c_fH}{2},\quad n_{Bt} = \frac{e_rl_e}{2mv_t}.
\end{align}

Now, we have presented the system of equations that can be used to compute all transition rates and vehicle numbers as functions of the following independent decision variables: the numbers of idle e-scooters at stations $\{n_{bs}\}$ and at random locations $\{n_{br}\}_{b > 0}$, station number $K^2$, station capacity $Q$, repositioning headway $H$ and initial truck load $R$, charging promotions $\{\pi_b\}$, and priority weights $\{\theta_{\hat{b}, b}\}$. 

\subsection{System Performance and Optimal System Design}
We can evaluate the system performance, e.g., agency cost and rider travel time. 
The agency costs are composed of four parts: (i) charging station installation cost $\omega_1\cdot K^2 + \omega_2\cdot K^2Q$; (ii) e-scooter fleet operation cost $\gamma\cdot F$, where $F = 
\left(\sum^{B}_{b = 0} (n_{br} + n_{bs}) + \sum^{B}_{b = 1} (n_{bw} + n_{bu}) + n_{0t} + n_{0f} + n_{Bt} + n_{Bf}\right)$ is the operating fleet size
; (iii) repositioning cost $\kappa\cdot l_e/H$, where $l_e/H$ is the truck distance units traveled per time unit
; and (iv) promotions paid to customers for dropping off e-scooters at charging stations $\sum_b\sum_{\hat{b}}\pi_b\cdot d_{(b+\hat{b})b,s}$.\footnote{Note that cost (iv) can also be interpreted as part of the riders' travel cost associated with the extra walking from charging stations to their destinations. Accordingly, we do not include this extra walking in the riders' cost.
} 
%
On the other hand, the average number of riders in the system is $\sum^{B}_{b = 1} (n_{bw} + n_{bu})$. Thus, the average rider travel time is $\sum^{B}_{b = 1} (n_{bw} + n_{bu})/(\lambda \Phi^2)$ based on Little's Law. Then, the average system-wide cost per rider, denoted by $Z$ [\$/rider], can be calculated as follows,
\begin{equation}\label{scooter_eq: system-wide cost}
    Z = \frac{1}{\lambda\Phi^2}\left[\omega_1\cdot K^2 + \omega_2\cdot K^2Q + \gamma\cdot F + \kappa\cdot \frac{l_e}{H}  + 
    \sum_{b}\sum_{\hat{b}}\pi_b\cdot d_{(b+\hat{b})b,s}+ \beta\sum^{B}_{b = 1} (n_{bw} + n_{bu})\right].
\end{equation}

The system design problem can be written as 
the following nonlinear programming model,
\begin{subequations}\label{scooter_eq:optimization}
\renewcommand{\theequation}{\theparentequation.\arabic{equation}}
\begin{align}
\min
\quad & Z
 \nonumber\\
\textrm{s.t.} \quad 
& \text{Eqn. } \eqref{scooter_eq: demand and trip length}-\eqref{scooter_eq: system-wide cost}, \nonumber \\
& a_{b\hat{b},r},\ a_{b\hat{b},s},\ p_b,\ d_{b(b-\hat{b}), r},\ d_{b(b-\hat{b}), s},\ c_b,\ c_f,\ e_i,\
n_{bi}\geq 0, \quad &&\forall b,\ \hat{b},\ i  \label{scooter_eq:constraint non-neg}\\
& 0\leq \pi_b \leq \beta S/v_w,\quad &&\forall b \in \{0, 1, ..., B-1\} \label{scooter_eq:constraint pi range}\\
& N_s \leq K^2Q, \label{scooter_eq:constraint capacity}\\
& \Phi/S_{\max} \leq K \leq \Phi/S_{\min}, \label{scooter_eq:constraint S range}\\
& Q_{\min} \leq Q \leq Q_{\max}, \\
& H_{\min} \leq H \leq H_{\max}, \\
& R_{\min} \leq R \leq R_{\max}, \label{scooter_eq:constraint R range} \\
& K,\ Q \in\mathbb{Z}^+,\label{scooter_eq:constraint K domain}\\
& \{n_{br}\},\ \{n_{bs}\},\ H,\ R,\ \{\pi_b\},\ \{\theta_{\hat{b}, b}\} \in \mathbb{R}^+.\label{scooter_eq:constraint other domain}
\end{align}
\end{subequations}
Constraints \eqref{scooter_eq:constraint non-neg} enforce all flow rates and fleet sizes to be non-negative real numbers; Constraints \eqref{scooter_eq:constraint pi range} set the range of promotions such that riders will only consider returning e-scooters to the stations nearest to their destinations; Constraint \eqref{scooter_eq:constraint capacity} ensures that the total number of e-scooters at stations do not exceed the total capacity of charging stations; Constraints \eqref{scooter_eq:constraint S range}-\eqref{scooter_eq:constraint R range} impose the range of decisions variables; and Constraints \eqref{scooter_eq:constraint K domain}-\eqref{scooter_eq:constraint other domain} define the domain of decision variables.

This optimization problem is highly nonlinear, and involves both continuous and integer decision variables. We propose the following solution approach: enumerate the value of $K$ over a feasible range, and for each $K$, treat $Q$ as a continuous variable, and use a gradient-based search method (with multiple starting solutions) to optimize the values of all continuous decision variables. 
The factorials in Eqn. \eqref{scooter_eq: prob station closer uncondition} and Eqn. \eqref{scooter_eq:pickup dist} 
may require careful numerical treatment.

\section{Numerical Studies}\label{scooter_sec:numerical studies}
In this section, we first verify the model predictions with agent-based simulations. Then, we conduct sensitivity analysis on a series of hypothetical examples to cast insights into the optimal system design under a variety of scenarios, and in particular, to highlight the benefits of deploying charging stations for shared e-scooter services.

We consider 1 du = 1 km and 1 tu = 1hr, and 
assume that $B = 8$ and $L_{\max} = 3$ km. 
The default parameter values are summarized in Table \ref{scooter_tab:parameters}. We use $v_s = 15$ km/hr as the cruising speed of shared e-scooters, an average between $9.13$ km/hr (which may include off-scooter time) in \citet{noland2019trip}, and 24 km/hr for commuting trips \citep{UNAGI}. 
Trucks speed is $v_t = 20$ km/hr, approximated by the commercial speed of buses \citep{Dickens_2023} considering that they both have many frequent stops.
The operation cost of shared e-scooters 
varies significantly across platforms and cities (e.g., 0.45 \$/veh-hr according to \citet{Jakobsen_2020} and 2.16 \$/veh-hr according to \citet{Griswold_2019}).\footnote{\citet{Jakobsen_2020} states that the operation cost of e-scooters is \$135 per e-scooter per month, including \$100 to lease the e-scooter and \$35 operation cost for software, maintenance, charging, etc. Assuming there are 30 days per month and 10 operation hours per day, this leads to an operation cost of $135/30/10 = 0.45$ \$/scooter-hr. \citet{Griswold_2019} shows that a shared e-scooter costs \$360 with a lifespan of one month. The operation cost (including charging, maintenance, customer support, etc.) is \$2.75 per ride. With 
an average demand of 3.49 trips per e-scooter per day, the average cost per e-scooter per hour is approximately \$2.16.} In this study, we use 
$\gamma = 1$ \$/veh-hr by default. The operation cost of charging stations $\omega_1 = 0.3$ \$/station-hr is borrowed from the cost per bus stop per hour \citep{ouyang2014continuum}. The cost of chargers is assumed to be $\omega_2 = \omega_1/5 = 0.06$ \$/charger-hr. The cost of trucks $\kappa = 4$ \$/truck-km is approximated by that of a bus \citep{ouyang2014continuum}, including $2$ \$/truck-km for 
vehicle acquisition, depreciation, fuel, maintenance, etc., and $2$ \$/truck-km for the crew (based on $40$ \$ per truck hour divided by the commercial speed $v_t = 20$ km/hr). We assume that it takes $8$ hr to fully charge a type-$0$ e-scooter, i.e., $\sum_b \tau_b = 8$, 
under two types of charging technologies: (i) a linear charging function, with $\tau_0 = \tau_1 = ... = \tau_7 = 1$ hr; and (ii) a piece-wise linear charging function, where the charging rate reduces by half when SoC exceeds a threshold of 80\% battery capacity 
\citep{montoya2017electric,osorio2021optimal}, with $\tau_0 = \tau_1 = ... = \tau_5 = 0.83$ hr, $\tau_6 = 1.33$ hr, and $\tau_7 = 1.67$ hr. 

\begin{table}[t]
    \centering
    \begin{adjustbox}{width=1\textwidth}
    \begin{tabular}{cccccccccc}
    \Xhline{2\arrayrulewidth}
         \multirow{2}{*}{Parameter} & $v_w$ & $v_s$  & $v_t$ & $\gamma$  & $\beta$  & $\omega_1$  & $\omega_2$  & $\kappa$ & $\sum_b \tau_b$\\
         & km/hr & km/hr & km/hr & \$/veh-hr & \$/hr & \$/station-hr & \$/charger-hr & \$/truck-km & hr \\ \hline
         Values & 3 & 15 & 20 & 1 & 20 & 0.3 & 0.06 & 4 & 8\\
    \Xhline{2\arrayrulewidth}
    \end{tabular}
    \end{adjustbox}
    \caption{Values of System Parameters.}
    \label{scooter_tab:parameters}
\end{table}

\subsection{Agent-based Simulations}\label{scooter_sec:simulation}
We develop agent-based simulations to verify the accuracy of the proposed analytical model, by comparing the expected passenger travel time from the simulation with that predicted by the model under the same setting. 

\begin{figure}[t]
     \centering
     \begin{subfigure}[b]{0.8\textwidth}
         \centering
         \includegraphics[width=\textwidth]{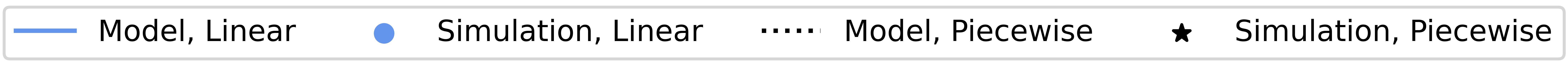}
     \end{subfigure}
     \begin{subfigure}[b]{0.49\textwidth}
         \centering
         \includegraphics[width=\textwidth]{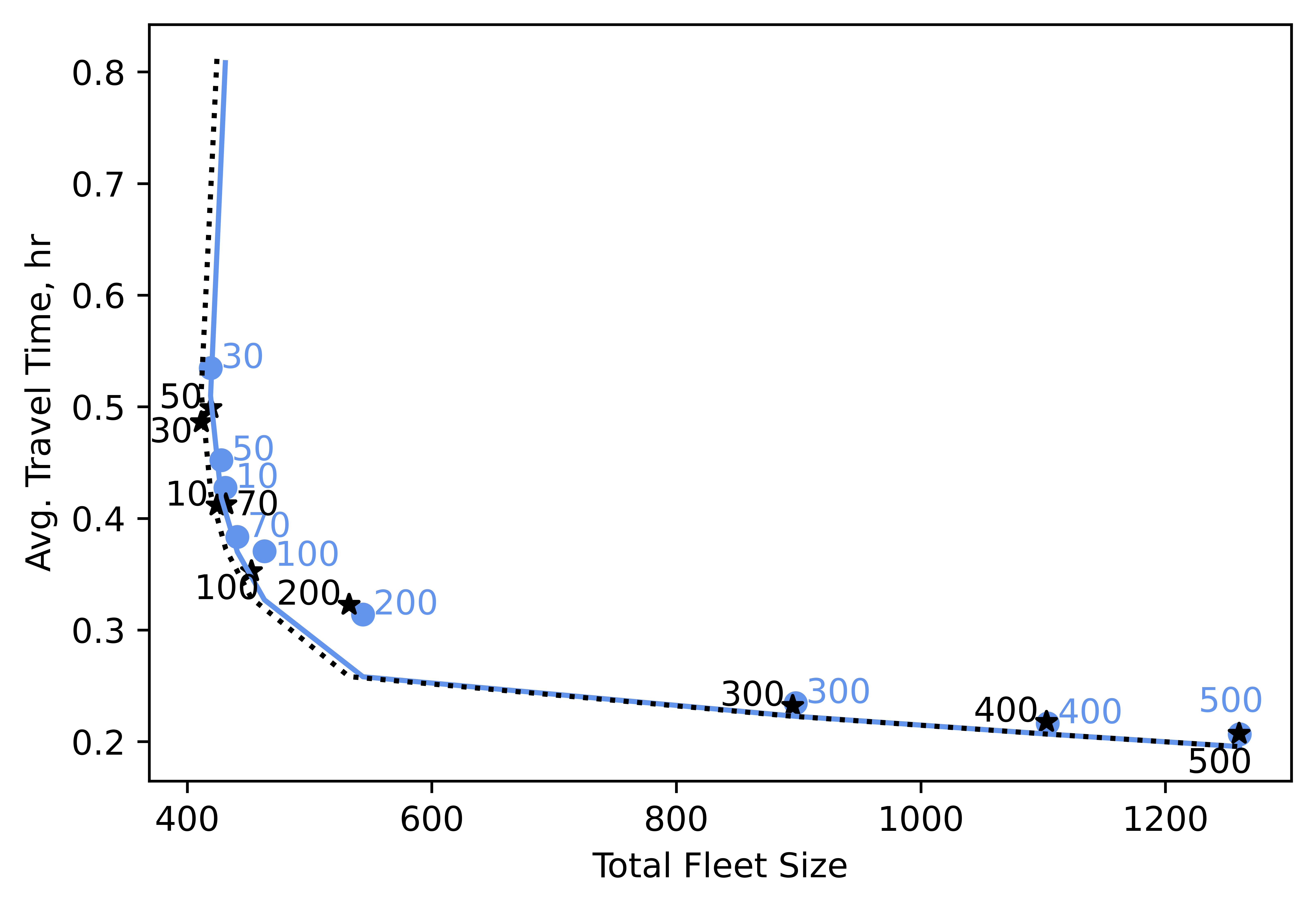}
         \caption{$K=5$.}
     \end{subfigure}
     \begin{subfigure}[b]{0.49\textwidth}
         \centering
         \includegraphics[width=\textwidth]{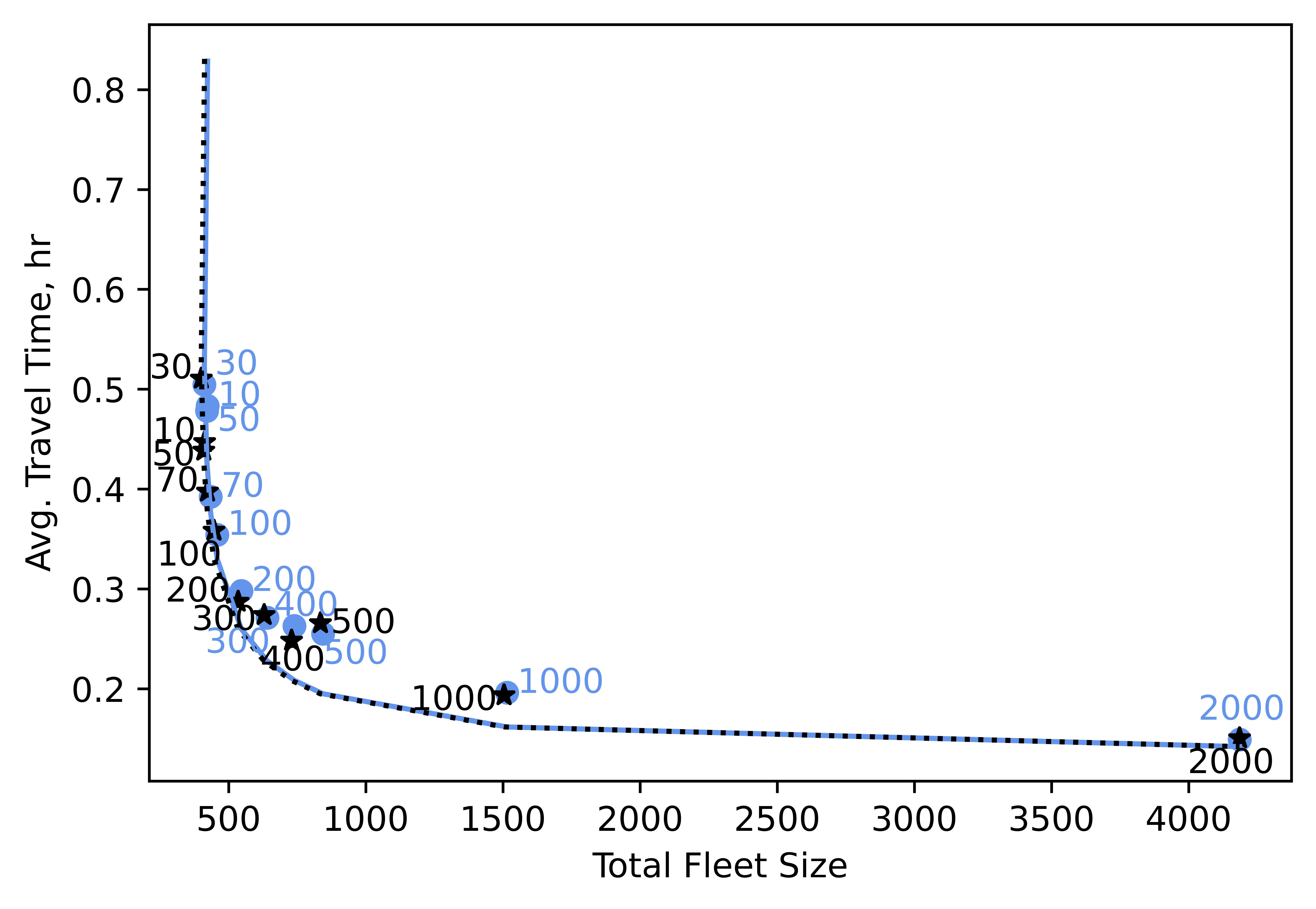}
         \caption{$K=10$.}
     \end{subfigure}
     \begin{subfigure}[b]{0.49\textwidth}
         \centering
         \includegraphics[width=\textwidth]{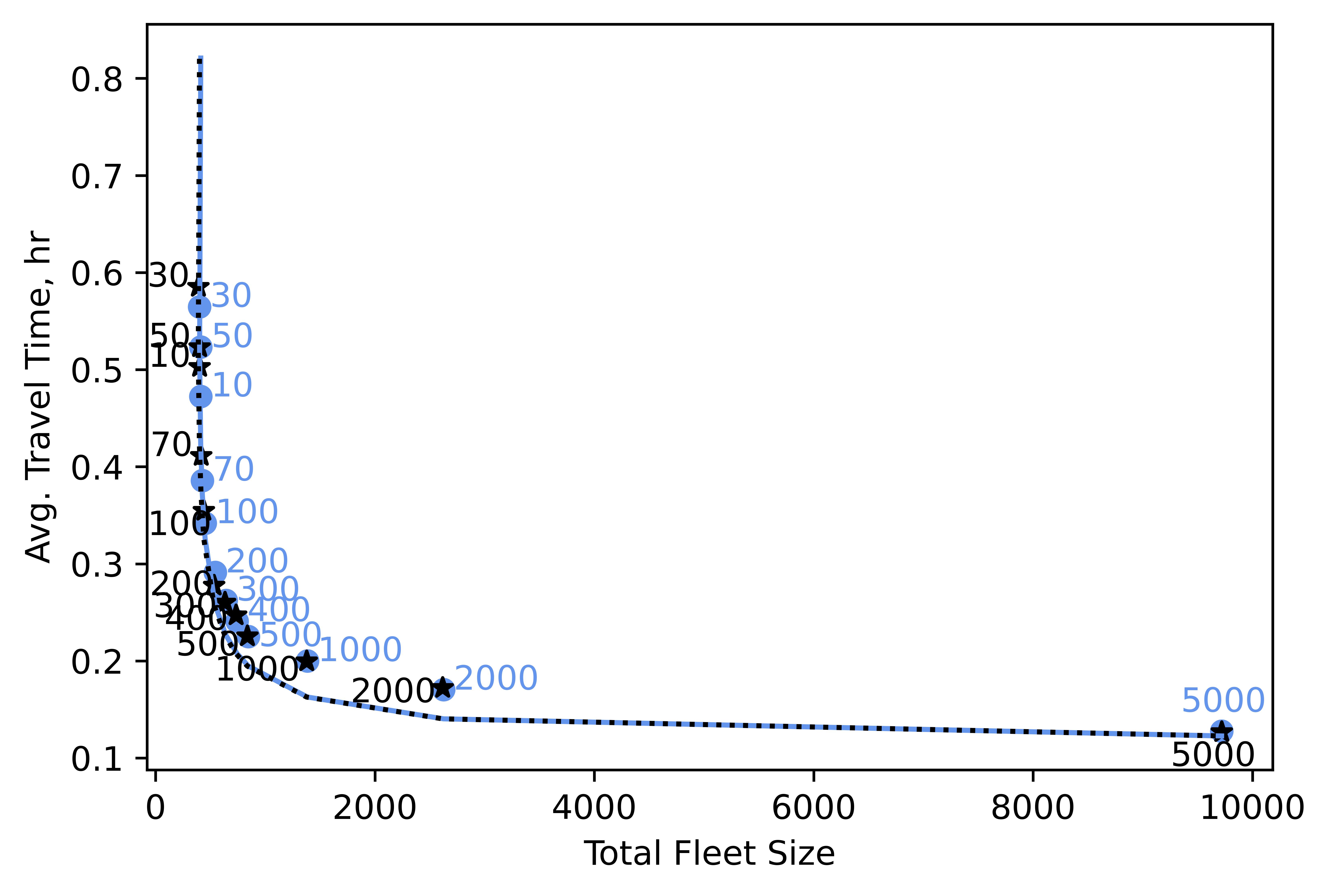}
         \caption{$K=15$.}
     \end{subfigure}
     \begin{subfigure}[b]{0.49\textwidth}
         \centering
         \includegraphics[width=\textwidth]{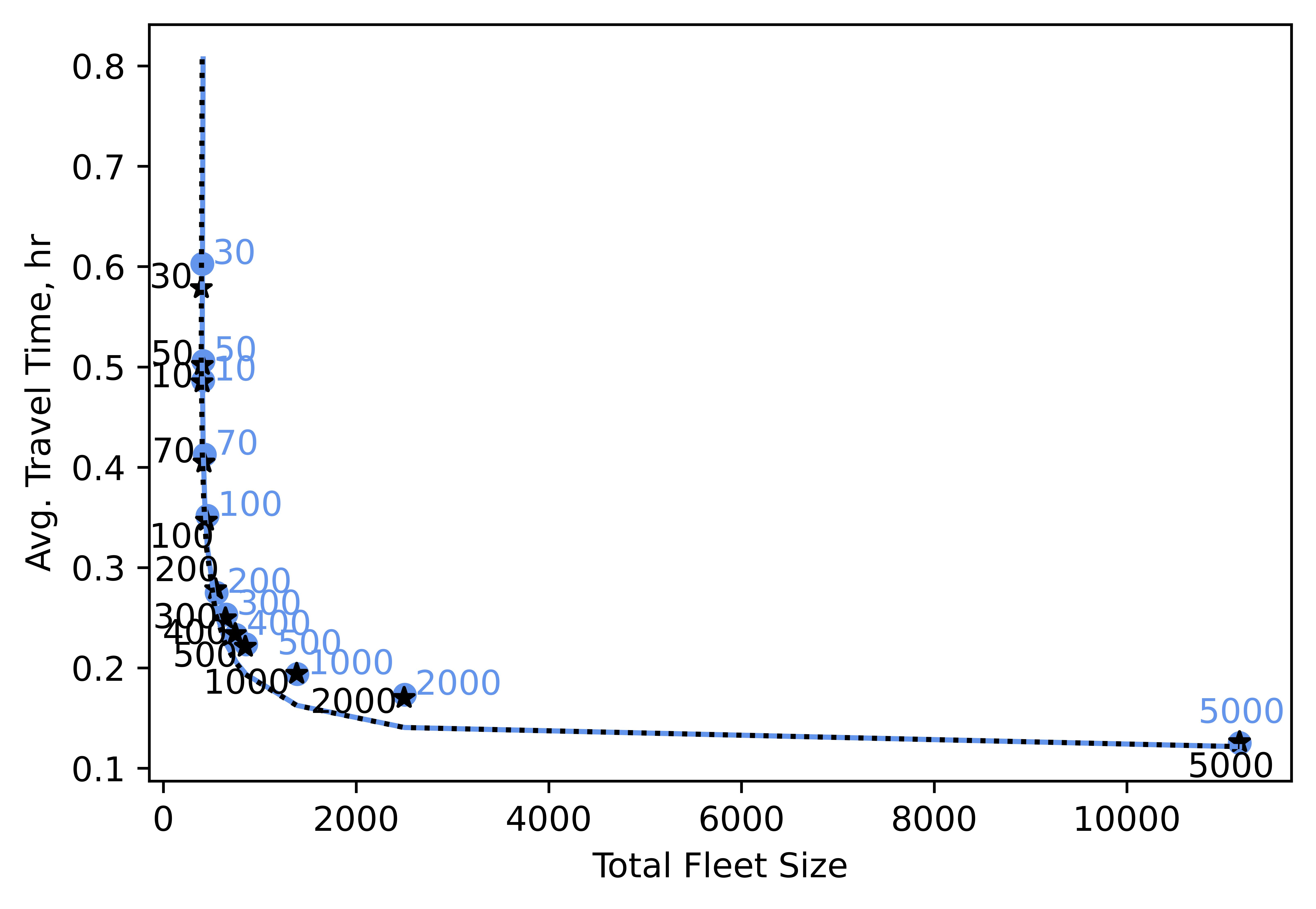}
         \caption{$K=20$.}
     \end{subfigure}
    \caption{Model predictions vs. simulation results.}
    \label{scooter_fig:SimResults}
\end{figure}

\begin{figure}[t]
     \centering
     \includegraphics[width=.5\textwidth]{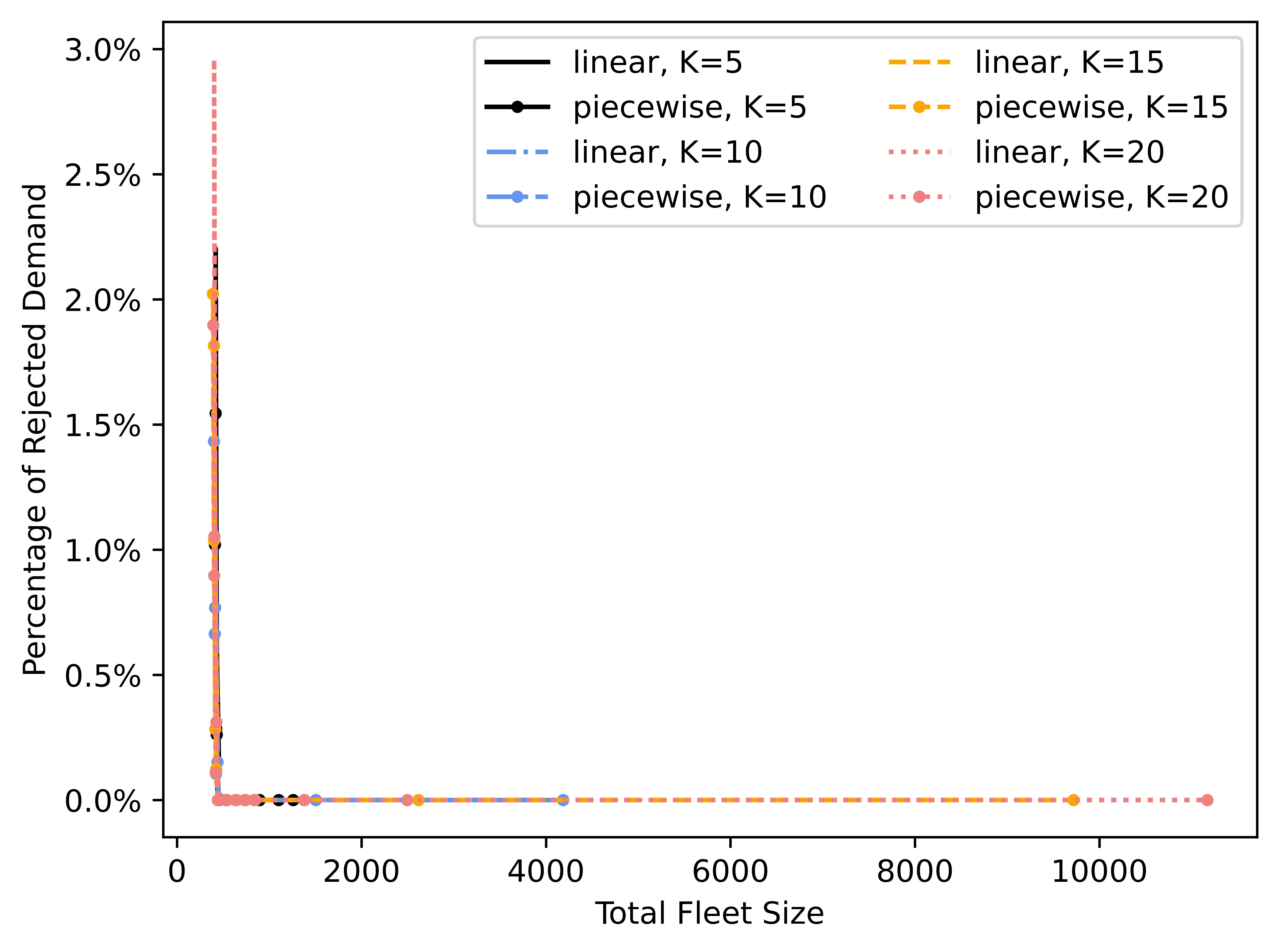}
    \caption{Proportion of rejected demand in simulation.}
    \label{scooter_fig:SimRejectedDemand}
\end{figure}

In addition to those in Table \ref{scooter_tab:parameters}, we set service region size $\Phi = 10$ km, and demand rate $\lambda = 1$ trip/hr-km$^2$. 
The depot is located $l_f = 2\Phi =20$ km south of the service region center. 
The value of $K$ varies from $5$ to $20$, and that of $\sum_{b>0} n_{br}$ from $10$ to $5000$, while all other decision variables are fixed: the truck headway $H = 1$ hr, truck load $R = 20$, the number of chargers per station $Q = 20$. The probability for the promotion of a type-$0$ and type-$1$ e-scooter to exceed the rider's walking disutility is fixed as $P_{\pi, 0} = 0.5$ and $P_{\pi, 1} = 0.25$, respectively, while no e-scooters at higher SoCs are dropped off at stations (i.e., $\pi_0 = \frac{\Phi\beta}{2Kv_w}$, $\pi_1 = \frac{\sqrt{2}\Phi\beta}{4Kv_w}$, and $\pi_b = 0,\ \forall b \geq 2$). In addition, the platform treats different types of e-scooters at charging stations indifferently as long as they are suitable for a rider, i.e., $\theta_{\hat{b}, b} = 1,\ \forall b\geq\hat{b}$. 

For each parameter combination, we first use the ``root" function in Scipy \citep{virtanen2020scipy} to solve the system of equations Eq. \eqref{scooter_eq: demand and trip length}-\eqref{scooter_eq:fleet size} and obtain the model predictions of the e-scooter quantity in each state and the expected passenger travel time. These e-scooter quantities are used to initialize the simulation. 
All e-scooters are initially idle with type-$B$ SoC, while $\sum_b n_{bs}$ e-scooters are randomly positioned among the $K^2$ stations and the rest randomly distributed within the service region. Trip requests are generated from a Poisson Process and served by the operating strategies in Section \ref{scooter_sec:model}. If no suitable e-scooters are available when a trip request arrives, we consider this trip request as lost. The truck routes are solved with a cluster-first route-second process. The clusters are obtained by solving a constrained k-means problem on the drop-off locations of type-$B$ e-scooters (which are randomly generated within the service area) to ensure the size of each cluster is less than the initial truck load $R = 20$, with the python package ``k\_means\_constrained" \citep{Levy-Kramer_k-means-constrained_2018}. The pickup locations of type-$0$ e-scooters are assigned to the clusters of their nearest centers (obtained from the abovementioned k-means algorithm). The route within each cluster is solved as a Traveling Salesman Problem using a simulated annealing algorithm implemented by the python package ``python-tsp" \citep{python-tsp}. 

Each simulation run simulates 2,000 hours of operations. The first 800 hours is the warm-up period and the last 200 hours is the cool-down period. Data from the remaining 1,000 hours, including the trip durations (from trip request to e-scooter drop-off) and the numbers of lost trips, are collected. The comparisons between model predictions and simulation results under linear and piece-wise linear charging functions are plotted in Figure \ref{scooter_fig:SimResults}, with the value of simulation input $\sum_{b>0} n_{br}$ labeled near the corresponding marker. The simulation results generally align very well with the model predictions, especially when the fleet size is large and the charging stations are dense. This makes sense, because the derivations of booking probabilities and pickup distances as in Eq. \eqref{scooter_eq:prob station closer by steps} and Eq. \eqref{scooter_eq:lp by steps} assume a large number of randomly-located idle vehicles -- i.e., the impacts of the boundary diminishes. 

Among the investigated parameter combinations, up to $3$\% of demand is lost when the fleet size is small, but this percentage quickly drops to zero as the fleet size increases. In addition, we found that the impacts of the linear vs. piece-wise linear charging functions 
are negligible, which indicates that the proposed model can predict system performance accurately under both types of charging functions. Therefore, in Section \ref{scooter_sec:sensitivity analysis}, we only use the piece-wise linear charging function in the sensitivity analysis.

In addition, the backward-bending curve of model predictions in Figure \ref{scooter_fig:SimResults}, especially for the case with $K=5$, indicates the existence of two equilibria under the same total fleet size --- the one on the lower branch is considered as a ``good" equilibrium, and the one on the upper branch is a ``bad" equilibrium, which is also known as the ``wild goose chase" phenomenon in taxi services \citep{castillo2017surge,Daganzo2019Paper}. The bad equilibria occur when the number of randomly located idle e-scooters is small, the expected pickup distance is large according to Eq. \eqref{scooter_eq:lp final}, and many e-scooters are trapped in the booked state waiting for riders to pick up. Meanwhile, it is interesting to observe that the backward-bending trend in model predictions 
becomes negligible when $K = 15$ and $K = 20$. This is reasonable, because the pickup distance is ``bounded" by the station spacing, and hence 
the wild goose chase phenomenon is less significant under smaller station spacings. This observation also partially explains why the model predictions deviate slightly more from simulations when $K$ is small --- as the simulated system states may include two types of equilibria, the measured system performance from simulation tends to be an affine combination of these two equilibria \citep{ouyang2023measurement}.

\subsection{Sensitivity Analysis}\label{scooter_sec:sensitivity analysis}
In this section, we apply the proposed analytical model to optimize the system design under the same system parameters (unless specified otherwise). We use solvers implemented in Scipy \citep{virtanen2020scipy}, Sequential Least Squares Programming algorithm \citep{kraft1988software} and Trust Region methods \citep{conn2000trust}, to obtain near-optimum solutions. The variable bounds are as follows: $S_{\min} = 0.5$ km, $S_{\max} = 5$ km, 
$Q_{\min} = 5$, $Q_{\max} = 20$, $H_{\min} = 1/6$ hr, $H_{\max} = 12$ hr, $R_{\min} = 5$, and $R_{\max} = 50$. For simplicity, we 
impose that $\pi_b = 0,\ \forall b \geq 4$, and only compare three types of priority weights (instead of directly optimizing them): (i) uniform weight (PW-1) for indifference across SoC, i.e., $\theta_{\hat{b}, b} = 1,\ \forall b \geq \hat{b}$; (ii) linear weight (PW-2) for higher SoC, i.e., $\theta_{\hat{b}, b} = b - \hat{b} + 1,\ \forall b \geq \hat{b}$; (iii) strongly nonlinear weight (PW-3) for near-full batteries,\footnote{PW-3 may be beneficial because it: (i) prohibits or discourages low-SoC e-scooters at stations to be picked up --- these e-scooters needs charging imminently, and (ii) 
encourages near-full e-scooters (i.e., $b\geq \floor{0.8B}$) to be picked up from the station --- these e-scooters may be charged at a lower rate than low-SoC ones. 
\label{scooter_fn:PW-3}} i.e., 
\begin{equation*}
    \theta_{\hat{b}, b} =\begin{cases}
        0,\quad & b = 1 \\
        b - \hat{b} + 1,\quad &\max\{2, \hat{b}\} \leq b < \floor{0.8B}\\
        10^{b - \hat{b} + 1},\quad &b \geq \floor{0.8B}
    \end{cases}.
\end{equation*}

The analysis focuses on optimizing the remaining decision variables: $\{n_{bs}\}$, $\{n_{br}\}_{b>0}$, $K$, $Q$, $\{\pi_b\}_{b<4}$, $H$, and $R$. We are interested in observing the optimal design characteristics, including the station spacing $S$, the density of idle e-scooters at random locations that are available for service $\sum_{b > 0} n_{br}/\Phi^2$, the average occupancy ratio of charging stations $\sum_b n_{bs}/(K^2Q)$, the average repositioning cost per trip $\frac{\kappa l_e}{\lambda \Phi^2 H}$, and the average charging incentive per trip $\sum_b\sum_{\hat{b}} \pi_b\cdot d_{(b+\hat{b})b, s}/(\lambda \Phi^2)$. In the meantime, we compare the system performance (e.g., agency cost per rider, average trip duration, and system-wide cost $Z$) of the proposed charging station design with that of two benchmark cases: a depot-only system (no charging stations) and a walk-only system (no e-scooters). The depot-only system is optimized with the proposed model, by removing the station-related states and transition rates in the queuing network in Figure \ref{scooter_fig:queuing network model}, and setting $K = Q = 0$ and $n_{bs} = \pi_b = 0,\ \forall b$.

\subsubsection{Varying Demand Rate}\label{scooter_sec:numerical_lambda}
\begin{table}
\begin{adjustbox}{width=\textwidth}
\begin{tabular}{ccccccccccccccc}
\Xhline{4\arrayrulewidth}
\multirow{2}{*}{$\lambda$} & \multicolumn{1}{c}{\multirow{2}{*}{\begin{tabular}[c]{@{}c@{}}Priority\\ Weights\end{tabular}}} & \multicolumn{7}{c}{Decision Variable} & \multicolumn{5}{c}{System Characteristics} \\ \cmidrule(lr){3-9}\cmidrule(lr){10-14}
& & $\sum\limits_{b > 0} n_{br}$ & $\sum\limits_b n_{bs}$ & $\{\pi_0,\ \pi_1,\ \pi_2,\ \pi_3\}$ & $K$ & $Q$ & $H$ & $R$ & $S$ & $\frac{\kappa l_e}{\lambda \Phi^2 H}$ & $\frac{\sum_b\sum_{\hat{b}} \pi_b\cdot d_{(b+\hat{b})b, s}}{\lambda \Phi^2}$ & $\frac{\sum_{b>0} n_{br}}{\Phi^2}$ & $\frac{\sum_b n_{bs}}{K^2Q}$ \\
\midrule
\multirow{3}{*}{1} & PW-1 & 91 & 292 & [4.96, 4.8, 3.48, 1.72] & 13 & 6 & 11.92 & 19 & 0.77 & 0.09 & 2.49 & 0.91 & 0.29 \\
 & PW-2 & 109 & 271 & [5.05, 4.82, 3.68, 1.93] & 13 & 7 & 11.37 & 30 & 0.77 & 0.03 & 2.24 & 1.09 & 0.23 \\
 & PW-3 & 124 & 249 & [4.69, 4.48, 3.82, 2.27] & 14 & 6 & 10.9 & 27 & 0.71 & 0.02 & 1.87 & 1.24 & 0.21 \\
\cline{1-14}
\multirow{3}{*}{5} & PW-1 & 278 & 1314 & [3.43, 3.26, 2.33, 1.21] & 19 & 11 & 11.66 & 18 & 0.53 & 0.02 & 1.98 & 2.78 & 0.33 \\
 & PW-2 & 321 & 1244 & [3.27, 3.14, 2.43, 1.45] & 20 & 10 & 11.64 & 16 & 0.5 & 0.02 & 1.72 & 3.21 & 0.31 \\
 & PW-3 & 420 & 1186 & [3.27, 3.15, 2.53, 1.83] & 20 & 9 & 11.99 & 18 & 0.5 & 0.02 & 1.5 & 4.2 & 0.33 \\
\cline{1-14}
\multirow{3}{*}{10} & PW-1 & 448 & 2573 & [3.2, 3.03, 1.99, 1.03] & 20 & 14 & 8.28 & 26 & 0.5 & 0.05 & 1.82 & 4.48 & 0.46 \\
 & PW-2 & 579 & 2514 & [3.18, 2.98, 2.11, 1.31] & 20 & 14 & 8.42 & 21 & 0.5 & 0.04 & 1.51 & 5.79 & 0.45 \\
 & PW-3 & 684 & 2369 & [3.22, 3.07, 2.2, 1.19] & 20 & 15 & 8.79 & 16 & 0.5 & 0.02 & 1.32 & 6.84 & 0.39 \\
\cline{1-14}
\multirow{3}{*}{50} & PW-1 & 2429 & 5563 & [1.66, 2.76, 1.35, 0.89] & 20 & 18 & 1.47 & 50 & 0.5 & 0.8 & 0.62 & 24.29 & 0.77 \\
 & PW-2 & 2878 & 5795 & [1.49, 2.9, 0.93, 0.96] & 20 & 19 & 1.49 & 50 & 0.5 & 0.77 & 0.53 & 28.78 & 0.76 \\
 & PW-3 & 2471 & 6225 & [1.57, 2.96, 1.39, 0.62] & 20 & 20 & 1.56 & 50 & 0.5 & 0.68 & 0.6 & 24.71 & 0.78 \\
\cline{1-14}
\multirow{3}{*}{100} & PW-1 & 4273 & 5260 & [0.9, 2.63, 0.0, 0.0] & 20 & 20 & 1.04 & 50 & 0.5 & 1.03 & 0.35 & 42.73 & 0.66 \\
 & PW-2 & 4441 & 5992 & [0.01, 2.73, 1.13, 0.59] & 20 & 20 & 1.05 & 50 & 0.5 & 0.99 & 0.34 & 44.41 & 0.75 \\
 & PW-3 & 4561 & 5808 & [0.04, 2.65, 1.06, 0.46] & 20 & 20 & 1.06 & 50 & 0.5 & 0.98 & 0.32 & 45.61 & 0.73 \\
\cline{1-14}
\bottomrule
\end{tabular}
\end{adjustbox}
\caption{Optimal design and system characteristics under varying demand rate $\lambda$.}
\label{scooter_tab: results when changing lambda}
\end{table}

\begin{figure}[t]
    \centering
    \includegraphics[width=\textwidth]{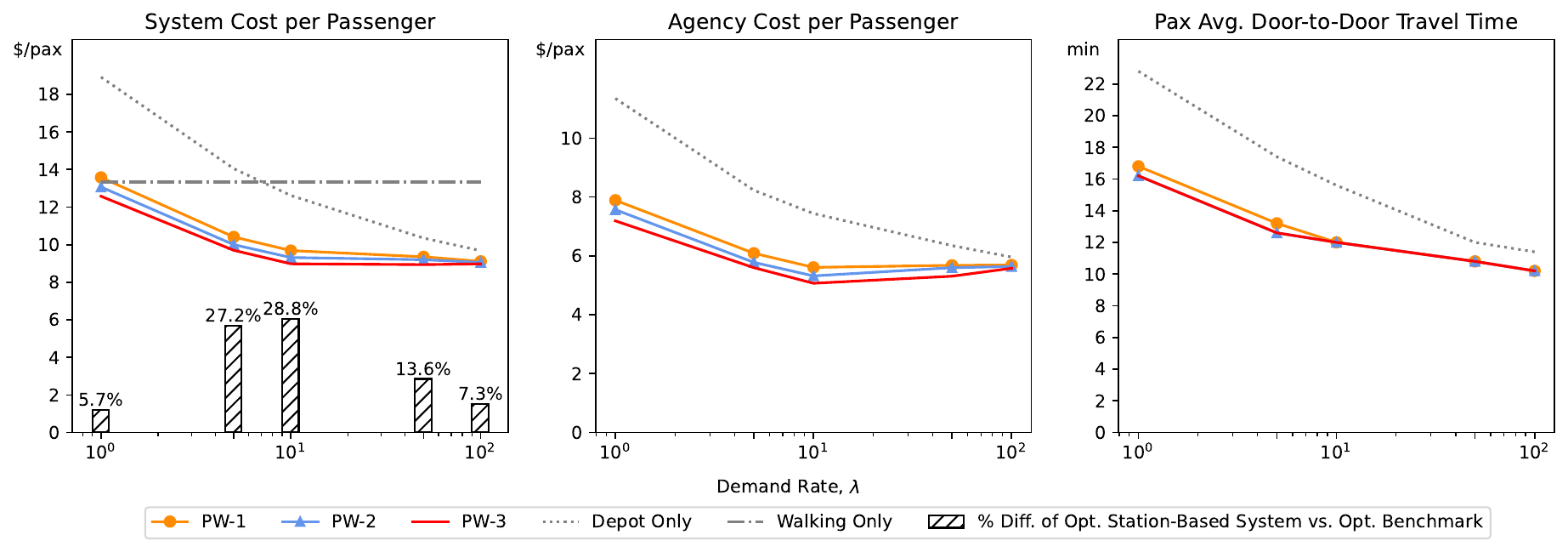}
    \caption{System performance under 
    varying demand rate $\lambda$.}
    \label{scooter_fig:varying lambda}
\end{figure}

In this subsection, we vary the demand rate $\lambda$ from $1$ to $100$ trip/hr-km$^2$. The optimal designs and system characteristics are summarized in Table \ref{scooter_tab: results when changing lambda}. Generally, as the demand rate increases, more idle e-scooters are deployed (indicated by increasing $\sum_{b > 0} n_{br}$ and $\sum_b n_{bs}$), the station density increases (indicated by greater $K$ and shorter station spacing $S$), and each station is equipped with more chargers $Q$. All these make sense, as larger demand can justify the operation cost of a larger pool of idle e-scooters and more charging infrastructure investment. More specifically, as demand increases, the ratio of randomly located idle e-scooters to those at stations $(\sum_{b > 0} n_{br})/(\sum_b n_{bs})$ first decreases when $\lambda \leq 10$ trip/hr-km$^2$ (i.e., more idle e-scooters are allocated to charging stations) and then increases (i.e., more idle e-scooters are allocated to random locations). It is possibly because that the total fleet size of idle e-scooters is small under a low demand rate, and the charging stations consolidate the e-scooters over space and provide riders with some ``guarantee" of access. When the demand rate is high, the total capacity at charging stations becomes limited but the total number of idle e-scooters grows very large, and therefore scattering them at random locations becomes beneficial. 

As demand increases, the optimal design also suggests reducing user-based charging operations (by lowering charging promotions $\{\pi_b\}$) but using more truck-based repositioning operations (e.g., shorter headway $H$ and larger initial truck load $R$). The major reasons could be that (i) truck-based repositioning becomes relatively more efficient due to spatial economies of scale, while providing the promotions imposes excessive financial burdens -- e.g., when $\lambda$ increases from 1 to 100 [trip/hr-km$^2$], the average promotion per trip $\sum_b\sum_{\hat{b}} \pi_b\cdot d_{(b+\hat{b})b, s}/(\lambda \Phi^2)$ is reduced by around \$2 while the truck repositioning cost per trip $\frac{\kappa l_e}{\lambda\Phi^2H}$ increases by less than \$1; and (ii) charging stations become congested under high demand rates (e.g., the average occupancy ratio of charging station $\sum_b n_{bs}/(K^2Q)$ increases from less than 0.3 when $\lambda = 1$ trip/hr-km$^2$ to over 0.7 when $\lambda \geq 50$ trip/hr-km$^2$), and therefore are not able to accommodate more e-scooters. 

In addition, it is interesting to observe the optimal $\{\pi_b\}$ across different types of e-scooters. When the demand rate $\lambda \leq 10$ trip/hr-km$^2$, optimal $\pi_b$ monotonically decreases with SoC $b$, i.e., $\pi_b > \pi_b'$ if $b < b'$, which is intuitive. However, when the demand rate $\lambda \geq 50$ trip/hr-km$^2$, the optimal designs apply a higher promotion of returning a type-$1$ e-scooter than that of a type-$0$ e-scooter. A possible reason is that, 
compared to a type-$0$ e-scooter, a type-$1$ e-scooter can be charged back to a suitable battery level faster, and it is beneficial to charge it at the station instead of at the depot. 

Across different priority weight types, we see that the optimal design with PW-3 tends to have more idle e-scooters at random locations and fewer idle e-scooters at stations, possibly due to the more efficient usage of e-scooters at charging stations under PW-3, as discussed in Footnote \ref{scooter_fn:PW-3}. As a result, the systems with PW-3 usually have the least average repositioning cost per trip
, the least average charging incentive payment per trip
, and the highest idle e-scooter density
. 

Figure \ref{scooter_fig:varying lambda} shows how the optimized e-scooter service with charging stations 
outperforms both benchmark cases under all investigated demand levels. We compare the optimal proposed system (among different priority weight types) with the better of the two benchmark systems, and the percentage differences in system-wide costs are represented by the bars in Figure \ref{scooter_fig:varying lambda}. First, among the proposed systems, the one with PW-3 always delivers the lowest system-wide cost, by saving $1.5\%\sim 7.4\%$ as compared to the one with PW-1 (which bears the highest system-wide cost). This difference indicates the potential benefits from also optimizing the priority weights in future studies. When the demand rate is extremely low ($\lambda = 1$ trip/hr-km$^2$), the depot-only system is obviously suboptimal, and the proposed system (with PW-3) is comparable with the walk-only system, with a slight saving of 5.7\% system-wide cost. This is reasonable, since the e-scooter trip lengths are suitable for walking, and there are little economies of scale from shared e-scooter services under extremely low demand. The system-wide costs of both the depot-only system and the proposed system drop as the demand rate increases, most notably for the depot-only system (due to the economies of scale from truck operations); as a result, when demand is very high ($\lambda = 100$ trip/hr-km$^2$), the gap between these two systems' costs almost diminishes. Under a moderate demand, the proposed system is significantly better than the two benchmarks, e.g., the proposed system (with PW-3) can save $27.2\%$ of the cost as compared to the walk-only system when $\lambda = 5$ trip/hr-km$^2$, and $28.8\%$ as compared to the depot-only system for the base case with $\lambda = 10$ trip/hr-km$^2$.  

In the following subsections, we fix the demand rate $\lambda = 10$ trip/hr-km$^2$ as a base case, and vary other system parameters one at a time for sensitivity analysis. 

\subsubsection{Varying Value-of-time}\label{scooter_sec:numerical_beta}

\begin{table}
\begin{adjustbox}{width=\textwidth}
\begin{tabular}{ccccccccccccccc}
\Xhline{4\arrayrulewidth}
\multirow{2}{*}{$\beta$} & \multicolumn{1}{c}{\multirow{2}{*}{\begin{tabular}[c]{@{}c@{}}Priority\\ Weights\end{tabular}}} & \multicolumn{7}{c}{Decision Variable} & \multicolumn{5}{c}{System Characteristics} \\ \cmidrule(lr){3-9}\cmidrule(lr){10-14}
& & $\sum\limits_{b > 0} n_{br}$ & $\sum\limits_b n_{bs}$ & $\{\pi_0,\ \pi_1,\ \pi_2,\ \pi_3\}$ & $K$ & $Q$ & $H$ & $R$ & $S$ & $\frac{\kappa l_e}{\lambda \Phi^2 H}$ & $\frac{\sum_b\sum_{\hat{b}} \pi_b\cdot d_{(b+\hat{b})b, s}}{\lambda \Phi^2}$ & $\frac{\sum_{b>0} n_{br}}{\Phi^2}$ & $\frac{\sum_b n_{bs}}{K^2Q}$ \\
\midrule
\multirow{3}{*}{10} & PW-1 & 279 & 2509 & [1.65, 1.57, 1.14, 0.62] & 20 & 15 & 11.56 & 21 & 0.5 & 0.02 & 1.06 & 2.79 & 0.42 \\
 & PW-2 & 325 & 2383 & [1.65, 1.58, 1.15, 0.62] & 20 & 14 & 11.61 & 22 & 0.5 & 0.02 & 0.99 & 3.25 & 0.43 \\
 & PW-3 & 509 & 2348 & [1.64, 1.59, 1.08, 0.85] & 20 & 14 & 11.29 & 20 & 0.5 & 0.01 & 0.77 & 5.09 & 0.42 \\
\cline{1-14}
\multirow{3}{*}{20} & PW-1 & 448 & 2573 & [3.2, 3.03, 1.99, 1.03] & 20 & 14 & 8.28 & 26 & 0.5 & 0.05 & 1.82 & 4.48 & 0.46 \\
 & PW-2 & 579 & 2514 & [3.18, 2.98, 2.11, 1.31] & 20 & 14 & 8.42 & 21 & 0.5 & 0.04 & 1.51 & 5.79 & 0.45 \\
 & PW-3 & 684 & 2369 & [3.22, 3.07, 2.2, 1.19] & 20 & 15 & 8.79 & 16 & 0.5 & 0.02 & 1.32 & 6.84 & 0.39 \\
\cline{1-14}
\multirow{3}{*}{30} & PW-1 & 786 & 2675 & [4.49, 4.08, 2.58, 0.64] & 20 & 14 & 6.59 & 36 & 0.5 & 0.1 & 2.19 & 7.86 & 0.48 \\
 & PW-2 & 731 & 2605 & [4.74, 4.46, 2.82, 1.47] & 20 & 15 & 11.59 & 44 & 0.5 & 0.02 & 2.04 & 7.31 & 0.43 \\
 & PW-3 & 923 & 2631 & [4.83, 4.48, 3.28, 2.53] & 20 & 15 & 11.42 & 17 & 0.5 & 0.01 & 1.55 & 9.23 & 0.44 \\
\cline{1-14}
\multirow{3}{*}{40} & PW-1 & 894 & 2836 & [5.24, 5.62, 3.24, 2.19] & 20 & 15 & 4.66 & 41 & 0.5 & 0.2 & 2.37 & 8.94 & 0.47 \\
 & PW-2 & 778 & 2588 & [5.76, 5.76, 4.03, 2.39] & 20 & 13 & 6.0 & 33 & 0.5 & 0.1 & 2.53 & 7.78 & 0.5 \\
 & PW-3 & 1027 & 2580 & [6.2, 5.72, 4.15, 2.56] & 20 & 14 & 10.0 & 33 & 0.5 & 0.03 & 1.99 & 10.27 & 0.46 \\
\cline{1-14}
\bottomrule
\end{tabular}
\end{adjustbox}
\caption{Optimal design and system characteristics under varying value-of-time $\beta$.}
\label{scooter_tab: results when changing beta}
\end{table}

\begin{figure}[t]
    \centering
    \includegraphics[width=\textwidth]{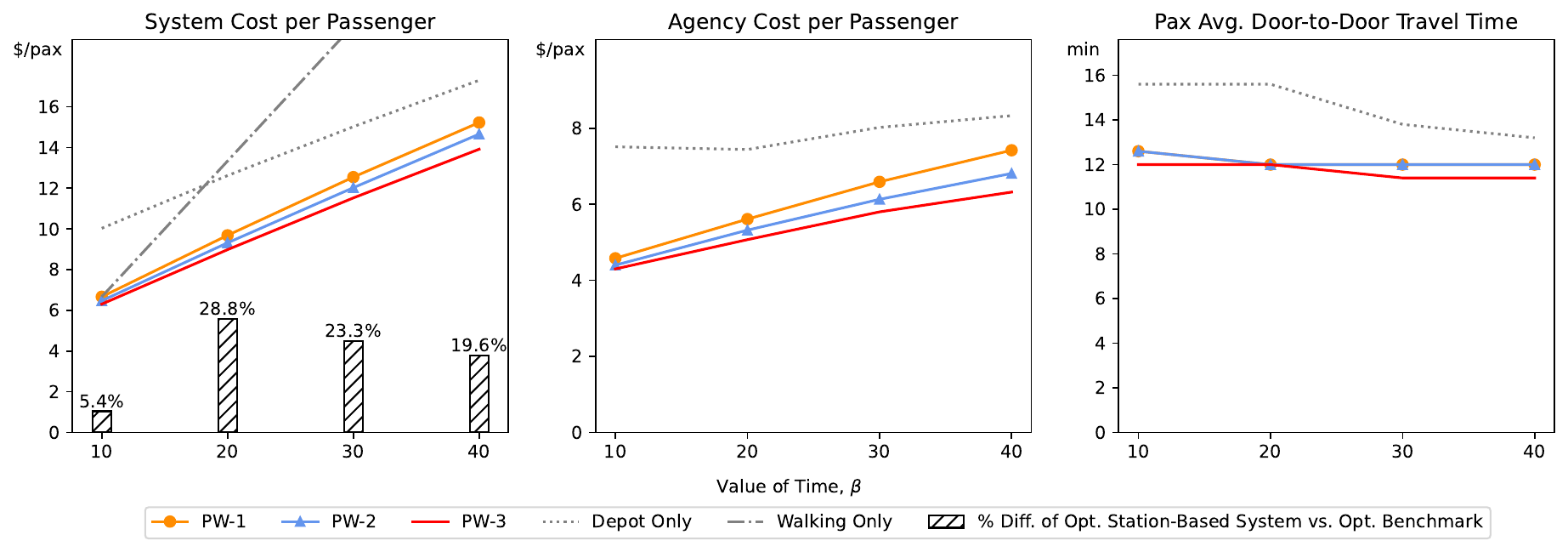}
    \caption{System performance under 
    varying value-of-time $\beta$.}
    \label{scooter_fig:varying beta}
\end{figure}

Table \ref{scooter_tab: results when changing beta} shows the optimal designs when the riders' value-of-time $\beta$ varies from $10$ \$/hr to $40$ \$/hr. As $\beta$ increases, more emphasis is put on rider travel experience, and hence more idle e-scooters are deployed at both random locations and at stations. This reduces the distance for a rider to pick up an e-scooter according to Eqn. \eqref{scooter_eq:lp final}. In addition, the repositioning headway $H$ tends to shrink so as to shorten the downtime of depleted e-scooters at random locations, at the expense of increasing repositioning cost $\frac{\kappa l_e}{\lambda\Phi^2 h}$. In the meantime, the number of charging stations $K^2$ and the number of chargers per station $Q$ are relatively stable, while the charging promotions $\{\pi_b\}$ are raised almost proportionally to $\beta$ so as to compensate for the increasing disutility of walking. 

Figure \ref{scooter_fig:varying beta} illustrates the comparison of system performance. As $\beta$ increases, both the depot-only system and the proposed system observe decreasing average rider travel time but more increases in agency cost, as expected. 
When $\beta = 10$ \$/hr, the proposed system (with PW-3) saves 5.4\% system-wide cost compared to the walk-only system. When $\beta \geq 20$ \$/hr, the advantage of the proposed system (with PW-3) is more significant, by reducing $19.6\%\sim 28.8\%$ system-wide cost as compared to the depot-only system. 


\subsubsection{Varying Cost Coefficients}\label{scooter_sec:numerical_cost}

In this subsection, we vary, one at a time, the operating cost of a shared e-scooter $\gamma$, the installation cost of a charging station $\omega_1$, and the repositioning cost $\kappa$. The optimal designs and system-wide characteristics are summarized in Tables \ref{scooter_tab: results when changing gamma}, \ref{scooter_tab: results when changing omega}, and \ref{scooter_tab: results when changing kappa}, respectively.

\begin{table}
\begin{adjustbox}{width=\textwidth}
\begin{tabular}{ccccccccccccccc}
\Xhline{4\arrayrulewidth}
\multirow{2}{*}{$\gamma$} & \multicolumn{1}{c}{\multirow{2}{*}{\begin{tabular}[c]{@{}c@{}}Priority\\ Weights\end{tabular}}} & \multicolumn{7}{c}{Decision Variable} & \multicolumn{5}{c}{System Characteristics} \\ \cmidrule(lr){3-9}\cmidrule(lr){10-14}
& & $\sum\limits_{b > 0} n_{br}$ & $\sum\limits_b n_{bs}$ & $\{\pi_0,\ \pi_1,\ \pi_2,\ \pi_3\}$ & $K$ & $Q$ & $H$ & $R$ & $S$ & $\frac{\kappa l_e}{\lambda \Phi^2 H}$ & $\frac{\sum_b\sum_{\hat{b}} \pi_b\cdot d_{(b+\hat{b})b, s}}{\lambda \Phi^2}$ & $\frac{\sum_{b>0} n_{br}}{\Phi^2}$ & $\frac{\sum_b n_{bs}}{K^2Q}$ \\
\midrule
\multirow{3}{*}{0.5} & PW-1 & 652 & 2733 & [3.11, 2.81, 1.83, 1.08] & 20 & 15 & 11.96 & 45 & 0.5 & 0.04 & 1.55 & 6.52 & 0.46 \\
 & PW-2 & 866 & 2788 & [3.19, 2.7, 1.91, 1.38] & 20 & 16 & 11.58 & 29 & 0.5 & 0.02 & 1.24 & 8.66 & 0.44 \\
 & PW-3 & 998 & 2644 & [3.16, 2.93, 2.14, 1.69] & 20 & 13 & 11.97 & 40 & 0.5 & 0.03 & 0.99 & 9.98 & 0.51 \\
\cline{1-14}
\multirow{3}{*}{1.0} & PW-1 & 448 & 2573 & [3.2, 3.03, 1.99, 1.03] & 20 & 14 & 8.28 & 26 & 0.5 & 0.05 & 1.82 & 4.48 & 0.46 \\
 & PW-2 & 579 & 2514 & [3.18, 2.98, 2.11, 1.31] & 20 & 14 & 8.42 & 21 & 0.5 & 0.04 & 1.51 & 5.79 & 0.45 \\
 & PW-3 & 684 & 2369 & [3.22, 3.07, 2.2, 1.19] & 20 & 15 & 8.79 & 16 & 0.5 & 0.02 & 1.32 & 6.84 & 0.39 \\
\cline{1-14}
\multirow{3}{*}{1.5} & PW-1 & 359 & 2530 & [3.28, 2.97, 2.09, 1.22] & 20 & 19 & 11.24 & 14 & 0.5 & 0.01 & 1.99 & 3.59 & 0.33 \\
 & PW-2 & 405 & 2459 & [3.27, 3.1, 2.34, 1.51] & 20 & 15 & 10.5 & 15 & 0.5 & 0.01 & 1.76 & 4.05 & 0.41 \\
 & PW-3 & 484 & 2348 & [3.26, 3.13, 2.44, 1.71] & 20 & 13 & 9.23 & 18 & 0.5 & 0.03 & 1.56 & 4.84 & 0.45 \\
\cline{1-14}
\bottomrule
\end{tabular}
\end{adjustbox}
\caption{Optimal design and system characteristics under varying e-scooter operation cost $\gamma$.}
\label{scooter_tab: results when changing gamma}
\end{table}

\begin{figure}[t]
    \centering
    \includegraphics[width=\textwidth]{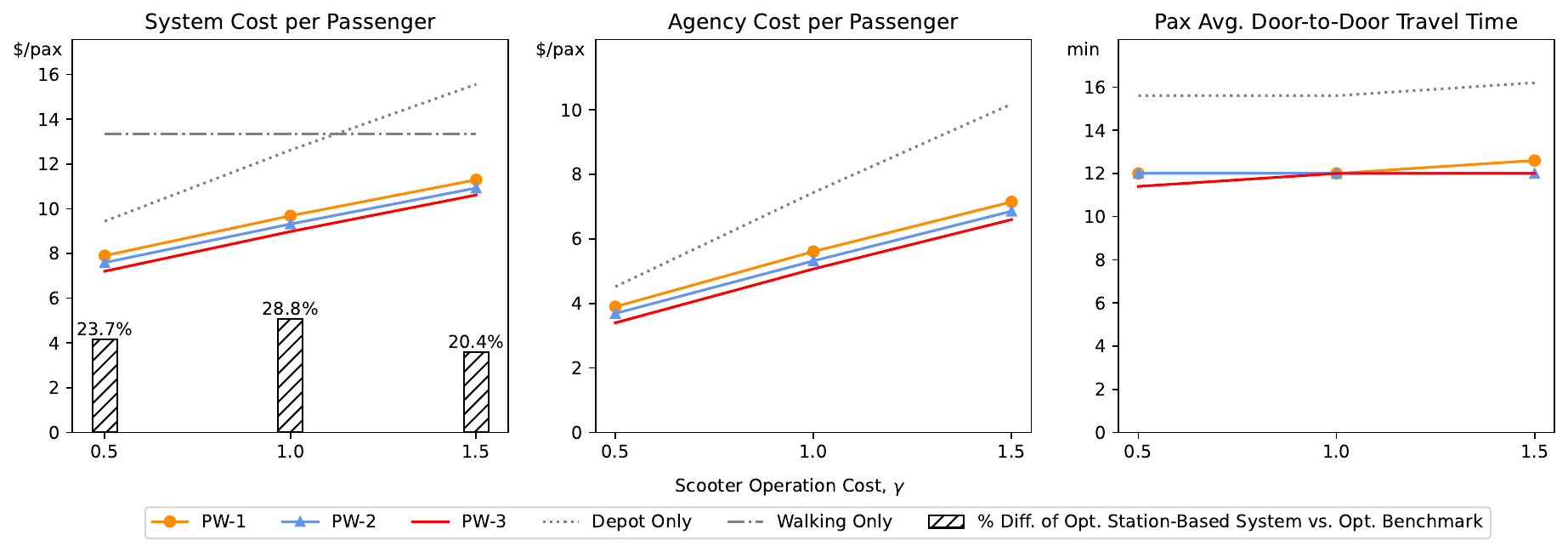}
    \caption{System performance under varying e-scooter operation cost $\gamma$.}
    \label{scooter_fig:varying gamma}
\end{figure}

As $\gamma$ increases, the optimal e-scooter fleet size tends to decrease
, and the system tends to expedite the charging process of e-scooters. On one hand, it offers higher promotions to encourage more drop-off of e-scooters at charging stations. On the other hand, it reduces the initial truck load $R$, which further reduces the number of e-scooters that are waiting for (or sitting on) the trucks.

Intuitively, one may expect headway $H$ to decrease with $\gamma$ too, which is consistent with the observation as $\beta$ goes from 0.5 \$/hr to 1.0 \$/hr. However, $H$ increases when $\gamma$ goes from 1.0 \$/hr to 1.5 \$/hr. A possible reason could be that more e-scooters are dropped off at stations when $\gamma = 1.5$ \$/hr, which balances (or even dwarfs) the impacts of extended headway $H$. As shown in Figure \ref{scooter_fig:varying gamma}, as $\gamma$ increases, both the agency costs and average trip duration of the depot-only system and the proposed system increase, which is reasonable. This figure also indicates that the depot-only system is very sensitive to $\beta$ while the proposed system is relatively insensitive. The reasons may include (i) riders can dynamically drop off and pick up e-scooters at stations, leading to more efficient utilization of e-scooters; and (ii) the proposed systems are less dependent on repositioning operations, which in turn reduces the number of e-scooters sitting at random locations, on trucks, and at the charging depot. By deploying the charging stations, the system-wide cost can be saved by 20.4\%$\sim$28.8\% as compared to the better of the two benchmark cases. 

\begin{table}
\begin{adjustbox}{width=\textwidth}
\begin{tabular}{ccccccccccccccc}
\Xhline{4\arrayrulewidth}
\multirow{2}{*}{$\omega_1$} & \multicolumn{1}{c}{\multirow{2}{*}{\begin{tabular}[c]{@{}c@{}}Priority\\ Weights\end{tabular}}} & \multicolumn{7}{c}{Decision Variable} & \multicolumn{5}{c}{System Characteristics} \\ \cmidrule(lr){3-9}\cmidrule(lr){10-14}
& & $\sum\limits_{b > 0} n_{br}$ & $\sum\limits_b n_{bs}$ & $\{\pi_0,\ \pi_1,\ \pi_2,\ \pi_3\}$ & $K$ & $Q$ & $H$ & $R$ & $S$ & $\frac{\kappa l_e}{\lambda \Phi^2 H}$ & $\frac{\sum_b\sum_{\hat{b}} \pi_b\cdot d_{(b+\hat{b})b, s}}{\lambda \Phi^2}$ & $\frac{\sum_{b>0} n_{br}}{\Phi^2}$ & $\frac{\sum_b n_{bs}}{K^2Q}$ \\
\midrule
\multirow{3}{*}{0.3} & PW-1 & 448 & 2573 & [3.2, 3.03, 1.99, 1.03] & 20 & 14 & 8.28 & 26 & 0.5 & 0.05 & 1.82 & 4.48 & 0.46 \\
 & PW-2 & 579 & 2514 & [3.18, 2.98, 2.11, 1.31] & 20 & 14 & 8.42 & 21 & 0.5 & 0.04 & 1.51 & 5.79 & 0.45 \\
 & PW-3 & 684 & 2369 & [3.22, 3.07, 2.2, 1.19] & 20 & 15 & 8.79 & 16 & 0.5 & 0.02 & 1.32 & 6.84 & 0.39 \\
\cline{1-14}
\multirow{3}{*}{0.6} & PW-1 & 516 & 2555 & [3.15, 2.99, 1.9, 0.9] & 20 & 13 & 6.88 & 32 & 0.5 & 0.08 & 1.75 & 5.16 & 0.49 \\
 & PW-2 & 709 & 2482 & [3.19, 2.92, 1.9, 0.8] & 20 & 13 & 8.58 & 35 & 0.5 & 0.04 & 1.45 & 7.09 & 0.48 \\
 & PW-3 & 877 & 2530 & [3.23, 3.1, 2.11, 1.56] & 20 & 13 & 11.84 & 28 & 0.5 & 0.03 & 1.13 & 8.77 & 0.49 \\
\cline{1-14}
\multirow{3}{*}{0.9} & PW-1 & 510 & 2488 & [3.1, 3.0, 1.93, 1.15] & 20 & 12 & 5.92 & 35 & 0.5 & 0.12 & 1.77 & 5.1 & 0.52 \\
 & PW-2 & 647 & 2450 & [3.2, 2.9, 1.97, 0.97] & 20 & 13 & 8.47 & 26 & 0.5 & 0.05 & 1.52 & 6.47 & 0.47 \\
 & PW-3 & 746 & 2568 & [3.6, 3.23, 2.42, 1.81] & 18 & 16 & 11.57 & 21 & 0.56 & 0.02 & 1.19 & 7.46 & 0.5 \\
\cline{1-14}
\bottomrule
\end{tabular}
\end{adjustbox}
\caption{Optimal design and system characteristics under varying charging station cost $\omega_1$.}
\label{scooter_tab: results when changing omega}
\end{table}

\begin{figure}[t]
    \centering
    \includegraphics[width=\textwidth]{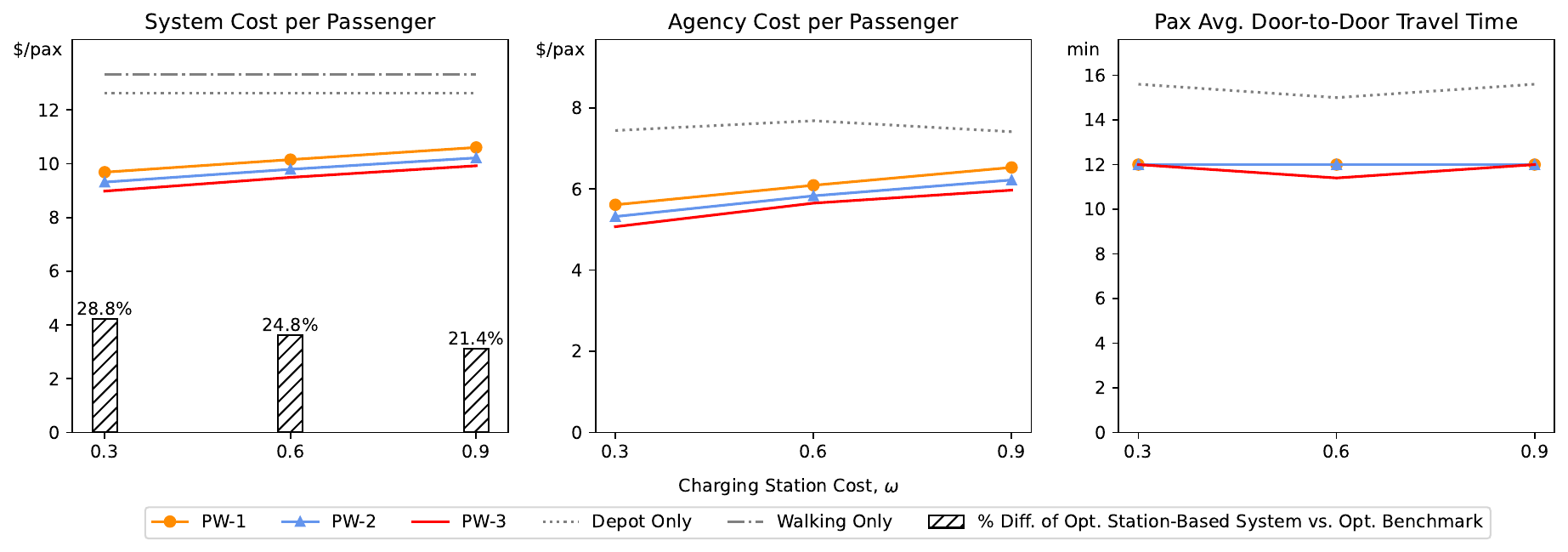}
    \caption{System performance under varying charging station cost $\omega_1$.}
    \label{scooter_fig:varying omega}
\end{figure}

The value of $\omega_1$ has almost negligible impacts on both the optimal designs and the system performances, as summarized in Table \ref{scooter_tab: results when changing omega} and illustrated in Figure \ref{scooter_fig:varying omega}. As $\omega_1$ increases, the optimal designs tend to slightly reduce the number of charging stations $K^2$ or the number of chargers per station $Q$, and the system-wide cost of station-based systems increases only moderately. This is possibly because the investment in charging stations is a relatively small portion of the system-wide cost, and hence the perturbation of its value does not drastically affect the system.

\begin{table}[t]
\begin{adjustbox}{width=\textwidth}
\begin{tabular}{ccccccccccccccc}
\Xhline{4\arrayrulewidth}
\multirow{2}{*}{$\kappa$} & \multicolumn{1}{c}{\multirow{2}{*}{\begin{tabular}[c]{@{}c@{}}Priority\\ Weights\end{tabular}}} & \multicolumn{7}{c}{Decision Variable} & \multicolumn{5}{c}{System Characteristics} \\ \cmidrule(lr){3-9}\cmidrule(lr){10-14}
& & $\sum\limits_{b > 0} n_{br}$ & $\sum\limits_b n_{bs}$ & $\{\pi_0,\ \pi_1,\ \pi_2,\ \pi_3\}$ & $K$ & $Q$ & $H$ & $R$ & $S$ & $\frac{\kappa l_e}{\lambda \Phi^2 H}$ & $\frac{\sum_b\sum_{\hat{b}} \pi_b\cdot d_{(b+\hat{b})b, s}}{\lambda \Phi^2}$ & $\frac{\sum_{b>0} n_{br}}{\Phi^2}$ & $\frac{\sum_b n_{bs}}{K^2Q}$ \\
\midrule
\multirow{3}{*}{2} & PW-1 & 475 & 2543 & [3.14, 2.97, 1.81, 0.99] & 20 & 15 & 5.67 & 18 & 0.5 & 0.03 & 1.8 & 4.75 & 0.42 \\
 & PW-2 & 514 & 2472 & [3.15, 3.04, 2.22, 1.43] & 20 & 13 & 5.51 & 17 & 0.5 & 0.03 & 1.6 & 5.14 & 0.48 \\
 & PW-3 & 753 & 2497 & [3.18, 3.17, 2.33, 1.73] & 20 & 13 & 6.72 & 15 & 0.5 & 0.02 & 1.23 & 7.53 & 0.48 \\
\cline{1-14}
\multirow{3}{*}{3} & PW-1 & 497 & 2555 & [3.2, 2.97, 1.71, 1.09] & 20 & 15 & 10.24 & 20 & 0.5 & 0.03 & 1.8 & 4.97 & 0.43 \\
 & PW-2 & 555 & 2458 & [3.17, 2.99, 2.08, 1.35] & 20 & 13 & 7.79 & 24 & 0.5 & 0.04 & 1.57 & 5.55 & 0.47 \\
 & PW-3 & 818 & 2554 & [3.26, 3.02, 2.27, 1.72] & 20 & 15 & 11.68 & 15 & 0.5 & 0.01 & 1.14 & 8.18 & 0.43 \\
\cline{1-14}
\multirow{3}{*}{4} & PW-1 & 448 & 2573 & [3.2, 3.03, 1.99, 1.03] & 20 & 14 & 8.28 & 26 & 0.5 & 0.05 & 1.82 & 4.48 & 0.46 \\
 & PW-2 & 579 & 2514 & [3.18, 2.98, 2.11, 1.31] & 20 & 14 & 8.42 & 21 & 0.5 & 0.04 & 1.51 & 5.79 & 0.45 \\
 & PW-3 & 684 & 2369 & [3.22, 3.07, 2.2, 1.19] & 20 & 15 & 8.79 & 16 & 0.5 & 0.02 & 1.32 & 6.84 & 0.39 \\
\cline{1-14}
\bottomrule
\end{tabular}
\end{adjustbox}
\caption{Optimal design and system characteristics under varying repositioning cost $\kappa$.}
\label{scooter_tab: results when changing kappa}
\end{table}

\begin{figure}[t]
    \centering
    \includegraphics[width=\textwidth]{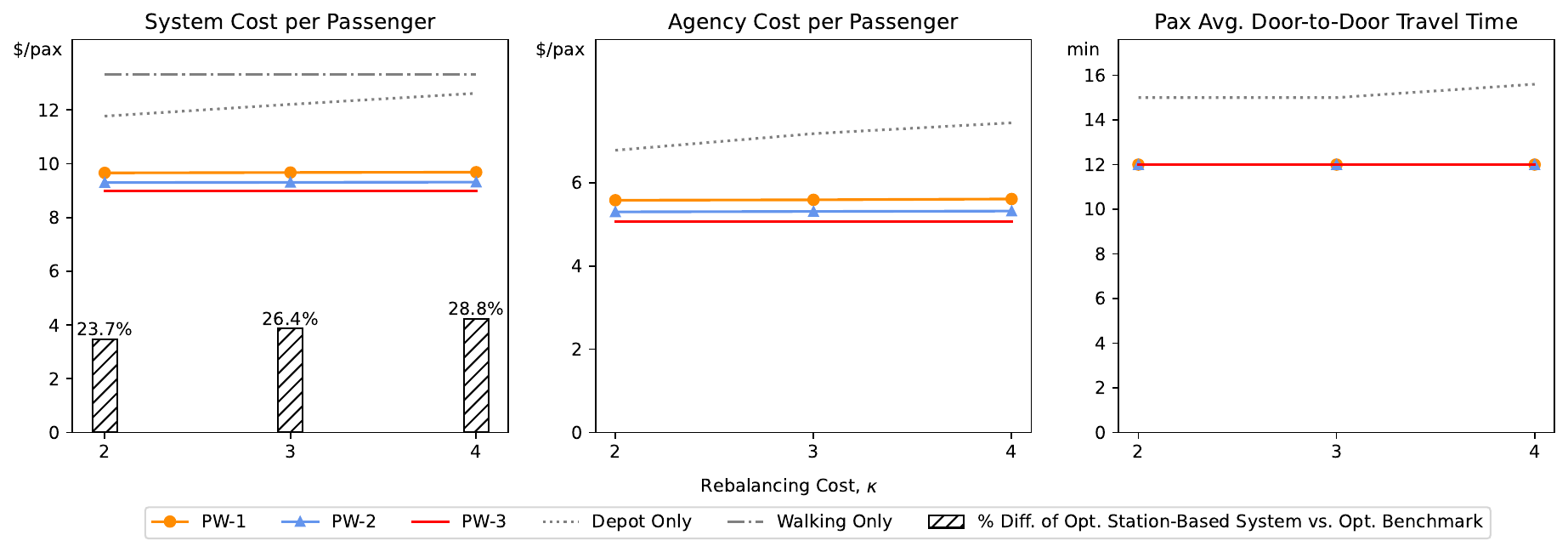}
    \caption{System performance under varying repositioning cost $\kappa$.}
    \label{scooter_fig:varying kappa}
\end{figure}

As $\kappa$ increases, the optimal headway $H$ and the optimal truck load $R$ tend to increase so as to reduce the total distance traveled per hour $l_e/H$, which is expected. In the meantime, the optimal charging promotion tends to slightly increase, so as to encourage more drop-offs at stations (and, in turn, to reduce truck-based repositioning operations). As a result, the system-wide costs of the proposed systems are relatively stable, while that of the depot-only system is moderately sensitive to $\kappa$; see Figure \ref{scooter_fig:varying kappa}.

\subsubsection{Varying E-scooter Battery Capacity}\label{scooter_sec:numerical_B}
\begin{table}[t]
\begin{adjustbox}{width=\textwidth}
\begin{tabular}{ccccccccccccccc}
\Xhline{4\arrayrulewidth}
\multirow{2}{*}{$B$} & \multicolumn{1}{c}{\multirow{2}{*}{\begin{tabular}[c]{@{}c@{}}Priority\\ Weights\end{tabular}}} & \multicolumn{7}{c}{Decision Variable} & \multicolumn{5}{c}{System Characteristics} \\ \cmidrule(lr){3-9}\cmidrule(lr){10-14}
& & $\sum\limits_{b > 0} n_{br}$ & $\sum\limits_b n_{bs}$ & $\{\pi_0,\ \pi_1,\ \pi_2,\ \pi_3\}$ & $K$ & $Q$ & $H$ & $R$ & $S$ & $\frac{\kappa l_e}{\lambda \Phi^2 H}$ & $\frac{\sum_b\sum_{\hat{b}} \pi_b\cdot d_{(b+\hat{b})b, s}}{\lambda \Phi^2}$ & $\frac{\sum_{b>0} n_{br}}{\Phi^2}$ & $\frac{\sum_b n_{bs}}{K^2Q}$ \\
\midrule
\multirow{3}{*}{8} & PW-1 & 448 & 2573 & [3.2, 3.03, 1.99, 1.03] & 20 & 14 & 8.28 & 26 & 0.5 & 0.05 & 1.82 & 4.48 & 0.46 \\
 & PW-2 & 579 & 2514 & [3.18, 2.98, 2.11, 1.31] & 20 & 14 & 8.43 & 21 & 0.5 & 0.03 & 1.51 & 5.79 & 0.45 \\
 & PW-3 & 684 & 2369 & [3.23, 3.07, 2.2, 1.19] & 20 & 15 & 9.09 & 16 & 0.5 & 0.01 & 1.33 & 6.84 & 0.39 \\
\cline{1-14}
\multirow{3}{*}{12} & PW-1 & 608 & 2405 & [3.2, 2.98, 2.23, 1.49] & 20 & 13 & 9.33 & 30 & 0.5 & 0.04 & 1.47 & 6.08 & 0.46 \\
 & PW-2 & 724 & 2308 & [3.27, 3.03, 2.33, 1.74] & 20 & 14 & 11.06 & 13 & 0.5 & 0.01 & 1.23 & 7.24 & 0.41 \\
 & PW-3 & 967 & 2306 & [3.62, 3.38, 2.7, 2.14] & 18 & 16 & 11.84 & 16 & 0.56 & 0.01 & 0.93 & 9.67 & 0.44 \\
\cline{1-14}
\multirow{3}{*}{16} & PW-1 & 806 & 2385 & [3.23, 3.11, 2.37, 1.72] & 20 & 14 & 8.07 & 19 & 0.5 & 0.02 & 1.17 & 8.06 & 0.43 \\
 & PW-2 & 893 & 2257 & [3.62, 3.4, 2.6, 2.08] & 18 & 16 & 11.34 & 11 & 0.56 & 0.01 & 0.99 & 8.93 & 0.44 \\
 & PW-3 & 1311 & 2394 & [4.08, 3.86, 3.15, 2.65] & 16 & 18 & 11.98 & 16 & 0.62 & 0.01 & 0.64 & 13.11 & 0.52 \\
\cline{1-14}
\bottomrule
\end{tabular}
\end{adjustbox}
\caption{Optimal design and system characteristics under varying e-scooter battery capacity $B$.}
\label{scooter_tab: results when changing B}
\end{table}

\begin{figure}[t]
    \centering
    \includegraphics[width=\textwidth]{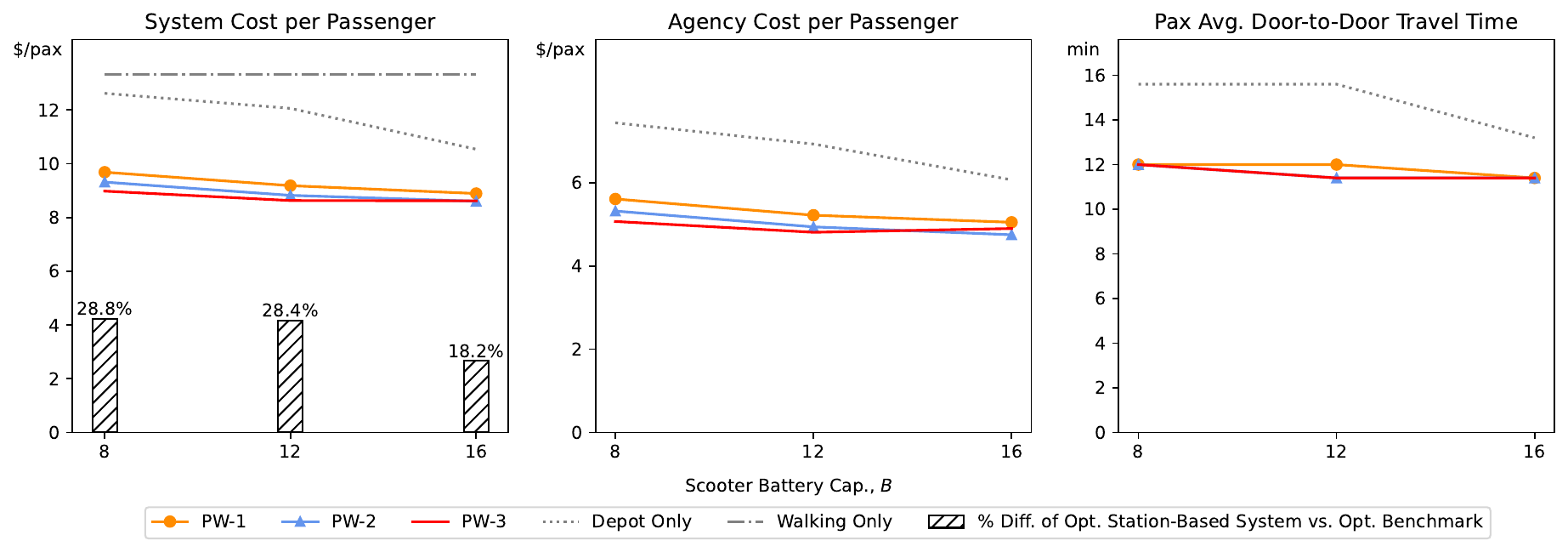}
    \caption{System performance under varying e-scooter battery capacity $B$.}
    \label{scooter_fig:varying B}
\end{figure}

In this subsection, we vary e-scooter battery capacity $B$, in terms of the travel range of a fully-charged e-scooter, from 8 km to 16 km, and extend the total charging time accordingly so that 
$\sum_b \tau_b = B$ hr. 
With larger battery capacities, e-scooters have less need for charging. As a result, the optimal number of charging stations $K^2$ tends to decrease and the optimal headway $H$ tends to increase, as summarized in Table \ref{scooter_tab: results when changing B}. As the station spacing $S$ increases, the optimal $\{\pi_b\}$ are increased such that the proportion of trips dropping off e-scooters at charging stations is relatively stable, and the optimal charger number $Q$ is expanded to accommodate more e-scooters at each station. Yet, with fewer charging stations, the total number of idle e-scooters at stations $\sum_b n_{bs}$ decreases, and idle e-scooters at random locations $\sum_{b > 0} n_{br}$ increases. Overall, both the average repositioning cost per trip and the average promotion cost per trip decrease as $B$ increases. 

The performance comparisons are plotted in Figure \ref{scooter_fig:varying B}. The benefit of a larger battery capacity on the system-wide cost is much less significant for the proposed system than for the depot-only system. This is expected. The expensive repositioning operations, which contribute to a major portion of the system-wide cost of the depot-only system, are directly impacted by the increase in battery capacity; e.g., the gap between the system-wide cost of the proposed system (with PW-3) and that of the depot-only system shrinks from 28.8\% to 18.2\% when the battery capacity increases from 8 km to 16 km.

\section{Conclusion}\label{scooter_sec:conclusion}
This paper proposes an aspatial queuing network model that jointly plans the services and charging stations (location and capacity) for a shared dockless electric micro-mobility system. To stay focused, the model is built with shared e-scooters as an example. Trips are categorized based on the required battery consumption, and riders are allowed to pick up and drop off e-scooters of different SoCs at both random locations and charging stations. The model incorporates two approaches to reposition and charge e-scooters: a user-based approach that incentivizes riders, with promotions, to drop off e-scooters at built charging stations, and a truck-based approach that collects and distributes e-scooters to and from an existing charging depot with trucks. E-scooters of different SoCs not only bear different priorities to be booked by the riders, but also are associated with different promotion values for drop-offs at charging stations. In the proposed queuing network model, an e-scooter's SoC and service status define its state, and the booking, pickup, drop-off, charging, and repositioning actions define state transitions. A system of nonlinear equations based on the queuing network is derived to describe the system's steady-state operations. The network model is further integrated into a constrained nonlinear programming model to optimize the system design, including e-scooter fleet sizes, station density and capacity, charging promotions, truck repositioning headway and truck load per dispatch, and priorities to utilize different types of e-scooters at stations. The proposed queuing network model is first verified with agent-based simulations, and then tested through a series of hypothetical scenarios to cast insights into the performance and optimal design of the proposed systems. The result indicates that the proposed system is superior to a walk-only system and a depot-only system in all investigated scenarios, especially when the customers have high value-of-time, the demand rate is moderate, e-scooters have relatively limited battery capacity with a moderate operating cost, and trucks are expensive to operate. As such, planning charging stations is a viable way to support the operations of dockless electric micro-mobility systems.

The proposed queuing network model for dockless e-scooter systems with charging stations can be further extended in various directions. We made a series of simplifying assumptions in this study so as to stay focused, but these assumptions can be relaxed. For example, we assumed that trucks only collect completely depleted e-scooters from, and distribute fully charged e-scooters to, random locations. Such limitations on the e-scooter types and locations can be easily removed by allowing trucks to collect (or distribute) e-scooters of more types of SoCs (e.g., lower or higher than any thresholds) and visit charging stations. The stations with truck visits usually have larger capacities to accommodate batch arrivals of e-scooters, and are often accessible by dedicated staff only, also known as ``charging hubs" or enterprise portals \citep{perch}. Therefore, an interesting extension is to adapt the proposed model for operations with charging hubs and possibly other charging approaches (e.g., battery swapping and use of independent contractors). 

The priority weights were limited to three specific cases in this study for simplicity, and the preliminary results demonstrated strong potential benefits from optimizing priority weights. Therefore, another extension could be done by numerically optimizing the priority weights (and implementing them via the apps) to further enhance the operation efficiency. Alternatively, the priority weights can also be explicitly calibrated with real-world data, to reveal riders' preferences for e-scooters at different SoCs in real-world case studies. 

In addition, this paper assumed that the demand is spatially homogeneous and hence balanced, which avoided the challenges of repositioning e-scooters between random locations. 
A possible solution to heterogeneous demand is to develop a zone-based model, similar to the one proposed in \citet{liu2023planning}, 
which explicitly determines the e-scooter flow within and between zones, rebalances e-scooters between zones, and optimizes the e-scooter deployment and station density within each zone. Furthermore, with a spatially heterogeneous demand, it would be interesting to explore location-dependent promotion mechanism design, which offers variable rewards to encourage e-scooter drop-offs at different stations for charging. 

Future studies can also look into many other dimensions, e.g., the riders' mode choices across two or more options (e.g., riding an e-scooter vs. walking vs. driving), the operations of other electric vehicles (e.g., e-bikes, drones, and shared electric cars), different charging technologies (e.g., fast charging vs. normal charging), temporally heterogeneous demand (e.g., peak hour vs. non-peak hour vs. non-operating hour), elastic demand and possibly third-party contractors (e.g., those who collect and reposition e-scooters to stations for a profit). 

\newpage
\appendix
\renewcommand*\thetable{\arabic{table}}
\renewcommand*\thefigure{\arabic{figure}}

\section{Key Notation}\label{scooter_sec:appendix notation}
\begin{longtable}{p{0.19\linewidth}p{0.79\linewidth}}
\Xhline{4\arrayrulewidth}
Notation & Description  \\ \hline
\\[-2ex]
\multicolumn{2}{l}{\textit{\textbf{Abbreviation, Indices, and Sets}}} \\
\\[-2ex]
du & Distance unit.\\
tu & Time unit.\\
SoC & State-of-Charge, measured in the number of du an e-scooter can travel in this paper.\\
$\mathcal{B}$ = \{0, 1, ..., B\} & Set of e-scooter SoCs, where $B$ is the battery capacity measured in du.\\
$b\in\mathcal{B}$ & Index of e-scooter type, based on the current SoC.\\
$\hat{b}\in\{1, ..., L_{\max}\}$ & Index of trip type, based on the trip length measured in du, where $L_{\max}$ [du] is the maximum travel range of shared e-scooter trips.\\
$i\in\{w, u, r, s, t, f\}$ & Index of e-scooter service status --- waiting for pickup, being used by a rider, idle at a random location, at a charging station, on a truck, or at the depot.\\
$C1$ & The condition that there are suitable e-scooters for a rider at its nearest station.\\
$C2$ & The condition that a station is closer to a rider's origin than any randomly located suitable e-scooter.\\
\hline
\\[-2ex]
\textit{\textbf{Parameters}} \\
\\[-2ex]
$\Phi$ [du] & Side length of the square service zone.\\
$v_s,\ v_w,\ v_t$ [du/tu] & E-scooter cruising speed, rider walking speed, and truck commercial speed.\\
$\gamma$ [\$/tu] & Operating cost of an e-scooter per time unit.\\
$\beta$ [\$/tu] & Customer's value-of-time.\\
$\omega_1,\ \omega_2$ [\$/tu] & Investments for each charging station and each charger amortized to each operational time unit, respectively.\\
$\kappa$ [\$/du] & Operation cost of a truck per distance unit traveled.\\
$\lambda,\ \lambda_{\hat{b}}$ [trip/tu-du$^2$] & Trip generation rates of all types of trips and of type-$\hat{b}$ trips, respectively.\\
$L_{\hat{b}}$ [du] & Average trip length of type-$\hat{b}$ trips.\\
$l_f$ [du] & Distance between the central charging depot and the region center.\\
$\tau_b$ [tu] & Time needed to charge an e-scooter at SoC $b$ to $b+1$.\\\hline
\\[-2ex]
\textit{\textbf{Variables}}\\
\\[-2ex]
$S$ [du] & Spacing between two adjacent charging stations.\\
$K\equiv \Phi/S$ & Ratio of $\Phi$ over $S$, and $K^2$ indicates the total number of charging stations.\\
$Q$ & Number of chargers installed at each charging station.\\
$\pi_b$ [\$/trip] & Promotions to encourage a rider to drop off an e-scooter with post-trip SoC $b$ to a station.\\
$H$ [tu] & Repositioning headway.\\
$R$ & Number of fully-charged e-scooters carried by each truck upon leaving the depot.\\
$\theta_{\hat{b}, b}$ & Priority weight for a type-$b$ e-scooter to be available to a type-$\hat{b}$ rider when there are multiple types of e-scooters equally distant from this rider's origin.\\
$n_{bi}$ & Number of type-$b$ e-scooters with service status $i$.\\
$a_{b\hat{b},r},a_{b\hat{b},s}$ [trip/tu] & Rates at which type-$b$ e-scooters at random locations and at stations are booked by type-$\hat{b}$ riders, respectively.\\
$a_{b,r},a_{b,s}$ [trip/tu] & Total booking rates of type-$b$ e-scooters at random locations and at stations by all rider types, respectively.\\
$p_b$ [trip/tu] & Rates for type-$b$ e-scooters to be picked up.\\
$d_{b(b-\hat{b}), r},\ d_{b(b-\hat{b}), s}$ [trip/tu] &  Rates for type-$b$ e-scooters to be dropped off after serving type-$\hat{b}$ trips, at charging stations and at random locations, respectively.\\
$c_b,\ c_f$ [trip/tu] &  Rate for type-$b$ e-scooters to be charged into type-$(b+1)$ at charging stations and for type-$0$ e-scooters to be charged into type-$B$ at the depot, respectively.\\
$e_f, e_r$ [trip/tu] &  Rates to reposition e-scooters from random locations to the depot, and from the depot to random locations, respectively.\\
$N_{\hat{b}s},\ N_{\hat{b}r}$ & Numbers of idle e-scooters at stations and those at random locations that are suitable for type-$\hat{b}$ trips, respectively.\\
$N_s$ & Total number of e-scooters at all stations.\\
$P_{C1, \hat{b}},\ P_{C2, \hat{b}}$  & Probabilitities of condition $C1$ and condition $C2$ 
for type-$\hat{b}$ riders, respectively.\\ 
$P_Q$ & Probability that a station has vacant chargers.\\
$P_{\pi, b}$ & Probability for the charging promotion of a type-$b$ e-scooter to exceed the rider's walking disutility.\\
$l$ [du] & General distance variable used in equation derivation, e.g., an e-scooter trip length, and distance from a rider's origin to the nearest station.\\
$l_e$ [du] & Total truck distance traveled to reposition e-scooters per headway.\\
$l_{\hat{b},p}$ [du] & Distance for a type-$\hat{b}$ rider to pick up an e-scooter, either at a station or a random location.\\
$m$ & Number of trucks dispatched per repositioning headway.\\
$F$ [tu/tu] & Total operating time unit of shared e-scooter fleet per tu.\\
$Z$ [\$/trip] & System-wide cost per trip.\\
\Xhline{4\arrayrulewidth}
\caption{List of key notation.}
\label{scooter_tab:notation}
\end{longtable}

\newpage
\nolinenumbers
\bibliographystyle{elsarticle-harv}\biboptions{authoryear}
\bibliography{references}

\end{document}